\documentclass[12pt]{amsart}

\usepackage[T1]{fontenc}
\usepackage{lmodern}

\usepackage[margin=1in]{geometry}

\usepackage{amssymb}
\usepackage[shortlabels]{enumitem}
\usepackage{thmtools, thm-restate}
\usepackage{graphicx,float}
\usepackage{subcaption}
\usepackage{tikz}
\usepackage{tikz-cd}
\usepackage{venndiagram}

\usetikzlibrary{
  backgrounds,
  decorations.pathmorphing,
  decorations.pathreplacing,
  shapes.arrows,
  arrows.meta,
  shapes.geometric,
  calc,
  positioning,
  3d
}

\usepackage{verbatim}
\usepackage{hyperref} 
\hypersetup{
  colorlinks=true,
  linkcolor=blue,   
  citecolor=blue,   
  urlcolor=blue,   
  filecolor=blue   
}
\usepackage[capitalise, noabbrev]{cleveref}

\newcommand{\Z}{\mathbb{Z}}
\newcommand{\R}{\mathbb{R}}
\newcommand{\K}{\mathbb{K}}

\newcommand{\Pp}{\mathcal{P}}
\newcommand{\st}{\colon}

\newcommand{\qand}{\quad \text{and} \quad}
\newcommand{\qfor}{\quad \text{for} \quad}

\DeclareMathOperator{\dist}{dist}
\DeclareMathOperator{\Ind}{Ind}
\DeclareMathOperator{\len}{len}
\DeclareMathOperator{\lk}{lk}
\DeclareMathOperator{\del}{del}
\DeclareMathOperator{\starr}{st}

\newtheorem{theorem}{Theorem}[section]

\newtheorem{proposition}[theorem]{Proposition}
\newtheorem{corollary}[theorem]{Corollary}
\newtheorem{lemma}[theorem]{Lemma}

\theoremstyle{definition}
\newtheorem{definition}[theorem]{Definition}

\newtheorem{remark}[theorem]{Remark}
\newtheorem{question}[theorem]{Question}

{\Alph{mainthmr}}

\title{Planar ternary graphs, flag spheres, and Delannoy polynomials}

\author[M. Bayer]{Margaret Bayer}
\address[M. Bayer]{Department of Mathematics, University of Kansas, Lawrence, KS, USA}
\email{bayer@ku.edu}

\author[R. Danner]{Richard Danner}
\address[R. Danner]{Department of Mathematics, University of Vermont, Burlington, VT, USA}
\email{rsdanner@uvm.edu}

\author[T. Holleben]{Thiago Holleben}
\address[T. Holleben]{Max Planck Institute for Mathematics in the Sciences, Leipzig, Germany}
\email{holleben@mis.mpg.de}

\author[M. Kramer]{Marie Kramer}
\address[M. Kramer]{Department of Mathematics \& Statistics,
Youngstown State University,
Youngstown, OH, USA}
\email{mkramer02@ysu.edu}

\author[Y. Yang]{Yirong Yang}
\address[Y. Yang]{School of Mathematics, Georgia Institute of Technology, Atlanta, GA, USA}
\email{yyang917@gatech.edu}

\begin{document}

\begin{abstract}
  In 2022 Kim showed that when a graph $G$ is ternary (without induced cycles of length divisible by three), its independence complex $\Ind(G)$ is either contractible or homotopy equivalent to a sphere. In this paper, we show that when $\Ind(G)$ is homotopy equivalent to a sphere of dimension $\dim \Ind(G)$, the complex is Gorenstein. Equivalently, $G$ is a $1$-well-covered graph. This answers a question by Faridi and Holleben.

We then focus on the independence complexes of Gorenstein planar ternary graphs. We prove that they are boundaries of vertex decomposable simplicial polytopes. We show that the transformations among these flag spheres using edge subdivisions and contractions can be modeled by the Hasse diagram of the partition refinement poset. In addition, their $h$-polynomials are products of Delannoy polynomials and thus real-rooted. Finally, we demonstrate a way to construct nonplanar Gorenstein ($1$-well-covered) ternary graphs from planar ones. 
\end{abstract}

\maketitle

\section{Introduction}

In this paper, we explore the algebraic, combinatorial, geometric, and
topological properties of the independence complexes of \emph{ternary} graphs,
graphs that do not contain induced cycles of length divisible by three.
Independence complexes of graphs are \textit{flag simplicial complexes}. 
One way to study the 
topological properties of these simplicial complexes is to use the translation into
algebraic properties of the associated Stanley--Reisner ring.  For example, 
Stanley \cite{Stanley} proved that the Stanley--Reisner ring of a homology sphere
is \textit{Gorenstein}.

It turns out that the Gorenstein property for independence complexes is
strongly related to properties of the graphs. Staples \cite{staples1979} introduced a class of graphs called \emph{1-well-covered}, denoted $W_2$; these are defined by a property of the independent sets of the graph. This class of graphs was further studied by Pinter \cite{pinter1995class}, and plays an important role in our work.
Hoang and Trung \cite{hoang2016characterization} proved that for
triangle-free graphs without isolated vertices, a graph $G$ is in $W_2$ if and 
only if the independence complex of $G$ is Gorenstein.
Trung \cite{trung2018characterization} determined the complete list of
connected planar graphs in $W_2$ (and hence with Gorenstein independence
complexes); the list contains exactly one graph for each independence number.
These graphs (with one exception) are in fact ternary.

Kim \cite{kim} proved that a graph is ternary if and only if
the independence complex of every induced subgraph is contractible or
homotopy equivalent to a sphere.  (This resolved a conjecture by
Engstr\"om \cite{Engstrom-KM}.)  This led Faridi and Holleben \cite{ternary}
to ask the question of the
connection between homotopy type and homology of these independence complexes.
For our purposes we phrase the question as follows.
\begin{question}\label{q:ternary}
    Let $G$ be a ternary graph without isolated vertices, $\Ind(G)$ its independence complex, and $d-1$ the dimension of $\Ind(G)$. Are the following two statements equivalent?
\begin{enumerate}
    \item $\Ind(G)$ is homotopy equivalent to a $(d-1)$-dimensional sphere.
    \item $\Ind(G)$ is a homology sphere.
\end{enumerate}
\end{question}

In this paper, we give an affirmative answer to this question. 

\medskip

\begin{restatable}{mainthm}{thmA}\label{thm:ternarygorenstein}

    Let $G$ be a ternary graph without isolated vertices. Then $\Ind(G)$ is a homology sphere if and only if $\Ind(G)$ is homotopy equivalent to $S^{d}$, where $d = \dim \Ind(G)$.

\end{restatable}

Moreover, using the characterization of planar graphs with independence complexes that are homology spheres (\cite{trung2018characterization}), we are able to provide further details on the independence complex of planar ternary graphs. 

In general, for $d \ge 4$, deciding if a simplicial $(d-1)$-sphere can be realized as the boundary of a $d$-polytope is a very hard problem. 
Our result below shows that being planar and ternary imposes not only strong topological restrictions on the independence complexes but also strong geometric ones.
\begin{restatable}{mainthm}{thmB}\label{t:onlygmintro}
    Let $G$ be a planar ternary graph such that $\Ind(G)$ is a homology sphere. Then $\Ind(G)$ is combinatorially equivalent to the boundary of a simplicial polytope and is vertex decomposable. 
\end{restatable}

\begin{remark}
Theorems \ref{thm:ternarygorenstein} and \ref{t:onlygmintro} together imply that if $G$ is a planar ternary graph whose independence complex is $(d-1)$-dimensional, then $\Ind(G)$ is homotopy equivalent to a $(d-1)$-dimensional sphere if and only if $\Ind(G)$ is the boundary of a vertex decomposable simplicial polytope. 
Moreover, the only nonternary planar graph $G$ for which $\Ind(G)$ is a homology sphere is the graph of the triangular prism, for which $\Ind(G)$ is the hexagon.
\end{remark}

We study how the independence complexes of planar ternary graphs are related by edge subdivisions.  Lutz and Nevo \cite{lutznevo16} showed that any two PL homeomorphic flag simplicial complexes can be connected by a sequence of edge subdivisions.  For $d$-dimensional flag PL spheres we thus have a connected graph we call $T_d$: the graph has a vertex for each flag PL $d$-sphere and an edge from one sphere to another that is obtained by an edge subdivision.

In this paper, we describe the restriction $H_d$ of $T_d$ to independence complexes of planar ternary graphs. In the statement below, let $\Pp_n$ denote the graph of the Hasse diagram of the partition refinement poset of the partitions of $n$.

\begin{restatable}{mainthm}{thmC}\label{thm:partitionrefine}
    The graph $H_{n-1}$ is isomorphic to $\Pp_n$. In particular, $\Pp_n$ is an induced subgraph of $T_{n-1}$ for every $n$.
\end{restatable}

We then focus on combinatorial properties of the flag spheres in~\cref{t:onlygmintro}. We begin by giving a combinatorial interpretation for their $h$-vectors, and as a consequence, using results from~\cite{WZC2019,M2019} we prove the following.

\begin{restatable}{mainthm}{thmD}\label{t:fhintro}
    
    Let $G$ be a planar graph such that $\Ind(G)$ is a homology sphere. Then the $h$-polynomial $h(\Ind(G), t)$ of $\Ind(G)$ is real-rooted. If additionally $G$ is ternary and connected, then $h(\Ind(G), t)$ is a Delannoy polynomial.
\end{restatable}

The paper is structured as follows.
In Section \ref{sec:Background} we recall important definitions and related results. In Section \ref{sec:PropertiesFlagSpheres}, we investigate ternary graphs and their independence complexes, proving Theorem \ref{thm:ternarygorenstein}. We introduce and delve into the properties of the graph family $\{G_m\}$, which gives all triangle-free connected planar Gorenstein graphs without isolated vertices. The section ends with a proof of Theorem \ref{t:onlygmintro}. In Section \ref{sec:FlipGraph}, we further examine the effect of edge subdivisions and contractions on the flag spheres arising from $G_m$ and prove \cref{thm:partitionrefine}. Section \ref{sec:VectorsPlanarFlagSpheres} is dedicated to proving Theorem \ref{t:fhintro}. Along the way, we give an explicit formula for the $h$-polynomial of $\Ind(G_m)$ and find a connection between the $h$-vector and the Delannoy numbers. We conclude this paper by posing some open questions and giving future research directions that might be of interest in Section \ref{sec:FutureWork}.

\section{Background}\label{sec:Background}

\subsection{Simplicial complexes} A \emph{simplicial complex} $\Delta$ on a vertex set $V$ is a collection of subsets of $V$ such that $\tau \subset \sigma \in \Delta$ implies $\tau \in \Delta$. Elements of $\Delta$ are called \emph{faces}, and maximal faces (with respect to inclusion) are called \emph{facets}. The number $\dim \sigma = |\sigma| - 1$ is called the \emph{dimension} of $\sigma$, and the dimension of $\Delta$ is the maximum dimension of its faces. If every facet of $\Delta$ has the same dimension, we say $\Delta$ is \emph{pure}. Moreover, $0$-dimensional faces are called vertices, and $1$-dimensional faces are called edges. For $0$- and $1$-dimensional faces, we often omit the usual set notation if it does not cause ambiguity. That is, if $\{x\}$ is a face of $\Delta$, we denote it by $x$, and if $\{x,y\}$ is a face of $\Delta$, we denote it by $xy$.

\subsection{Graphs and independence complexes} 
\begin{figure}[H]
    \centering
    \resizebox{1\textwidth}{!}{
\begin{tikzpicture}[scale=1.5, line cap=round, line join=round, every node/.style={transform shape, font=\large}]

\begin{scope}[shift={(-7,0)}]
\begin{scope}

\coordinate (ai-2) at (-4,0);
\coordinate (ai-1) at (-2,0);
\coordinate (ai) at (0,0);
\coordinate (ai+1) at (2,0);
\coordinate (ai+2) at (4,0);

\coordinate (bi-2) at (-4,2);
\coordinate (bi-1) at (-2,2);
\coordinate (bi) at (0,2);
\coordinate (bi+1) at (2,2);
\coordinate (bi+2) at (4,2);

\coordinate (ci) at (-1,3);
\coordinate (ci+1) at (1,3);
\coordinate (ci+2) at (3,3);


\draw[thick] (ai) -- (bi);
\draw[thick] (ai+1) -- (bi+1);
\draw[thick] (ai+2) -- (bi+2);

\draw[thick] (ai) -- (bi);
\draw[thick] (bi+1) -- (ci+2) -- (bi+2);

\draw[thick] (bi+1) -- (ci+2);
\draw[thick] (ai+1) -- (ai+2);


\node[fill=black, circle, inner sep=2pt, label=-90:{$a_1$}] at (ai) {};
\node[fill=black, circle, inner sep=2pt, label=-90:{$a_1$}] at (ai) {};
\node[fill=black, circle, inner sep=2pt, label=-90:{$a_2$}] at (ai+1) {};
\node[fill=black, circle, inner sep=2pt, label=-90:{$a_3$}] at (ai+2) {};

\node[fill=black, circle, inner sep=2pt, label=0:{$b_1$}] at (bi) {};
\node[fill=black, circle, inner sep=2pt, label=0:{$b_2$}] at (bi+1) {};
\node[fill=black, circle, inner sep=2pt, label=0:{$b_3$}] at (bi+2) {};

\node[fill=black, circle, inner sep=2pt, label=90:{$c_3$}] at (ci+2) {};

\end{scope}

\begin{scope}[shift={(9,1)}]

\coordinate (T) at (0,2.5);    
\coordinate (B) at (0,-2.5);   
\coordinate (L) at (-2,0);   
\coordinate (R) at (2,0);    
\coordinate (F) at (.3,-0.45);  
\coordinate (K) at (-.3,0.35); 
\coordinate (c2) at (-1,-1.3);  
\coordinate (c3) at (-.83,-.22); 

\draw[thick] (T) -- (L) -- (F) -- (R) -- (T);
\draw[thick] (L) -- (B) -- (R);
\draw[thick] (F) -- (B);
\draw[thick] (F) -- (T);
\draw[thick] (c3) -- (T);
\draw[thick] (c3) -- (B);

\draw[thick, dashed] (T) -- (K);
\draw[thick, dashed] (T) -- (F);
\draw[thick, dashed] (L) -- (K);
\draw[thick, dashed] (R) -- (K);
\draw[thick, dashed] (B) -- (K);

\node[fill=black, circle, inner sep=2pt, label=above:$b_1$] at (T) {};
\node[fill=black, circle, inner sep=2pt, label=below:$a_1$] at (B) {};
\node[fill=black, circle, inner sep=2pt, label=left:$a_2$]  at (L) {};
\node[fill=black, circle, inner sep=2pt, label=right:$b_2$] at (R) {};
\node[fill=black, circle, inner sep=2pt, label=-30:$a_3$] at (F) {};
\node[fill=black, circle, inner sep=2pt, label=60:$b_3$] at (K) {};
\node[fill=black, circle, inner sep=2pt, label=210:$c_3$] at (c3) {};

\end{scope}
\end{scope}

\begin{scope}[shift={(9,0)}]

    \coordinate (ai) at (0,0);
    \coordinate (ai+1) at (2,0);
    \coordinate (ai+2) at (4,0);

    \coordinate (bi) at (0,2);
    \coordinate (bi+1) at (2,2);
    \coordinate (bi+2) at (4,2);
    
    \coordinate (ci+1) at (1,3);
    \coordinate (ci+2) at (3,3);


    \draw[thick] (ai) -- (bi);
    \draw[thick] (ai+1) -- (bi+1);
    \draw[thick] (ai+2) -- (bi+2);

    \draw[thick] (ai) -- (bi) -- (ci+1) -- (bi+1) -- (ci+2) -- (bi+2);

    \draw[thick] (ci+1) .. controls (1,.5) and (3,1.3) .. (ai+2);
    \draw[thick] (ci+1) .. controls (1,.5) and (3,1.3) .. (ai+2);

    \draw[thick] (ci+1) -- (bi+1) -- (ci+2);
    \draw[thick] (ai) -- (ai+1) -- (ai+2);

\node[fill=black, circle, inner sep=2pt, label=-90:{$a_1$}] at (ai) {};
\node[fill=black, circle, inner sep=2pt, label=-90:{$a_2$}] at (ai+1) {};
\node[fill=black, circle, inner sep=2pt, label=-90:{$a_3$}] at (ai+2) {};

\node[fill=black, circle, inner sep=2pt, label=180:{$b_1$}] at (bi) {};
\node[fill=black, circle, inner sep=2pt, label=0:{$b_2$}] at (bi+1) {};
\node[fill=black, circle, inner sep=2pt, label=0:{$b_3$}] at (bi+2) {};

\node[fill=black, circle, inner sep=2pt, label=90:{$c_2$}] at (ci+1) {};
\node[fill=black, circle, inner sep=2pt, label=90:{$c_3$}] at (ci+2) {};

  \end{scope}

\begin{scope}[shift={(18,1)}]
\coordinate (T) at (0,2.5);    
\coordinate (B) at (0,-2.5);   
\coordinate (L) at (-2,0);   
\coordinate (R) at (2,0);    
\coordinate (F) at (.3,-0.45);  
\coordinate (K) at (-.3,0.35); 
\coordinate (c2) at (-1,-1.3);  
\coordinate (c3) at (-.83,-.22); 

\draw[thick] (T) -- (L) -- (F) -- (R) -- (T);
\draw[thick] (L) -- (B) -- (R);
\draw[thick] (F) -- (B);
\draw[thick] (F) -- (T);
\draw[thick] (c2) -- (c3);
\draw[thick] (c3) -- (T);
\draw[thick] (c3) -- (B);

\draw[thick, dashed] (T) -- (K);
\draw[thick, dashed] (T) -- (F);
\draw[thick, dashed] (L) -- (K);
\draw[thick, dashed] (R) -- (K);
\draw[thick, dashed] (B) -- (K);
\draw[thick, dashed] (c2) -- (K);

\node[fill=black, circle, inner sep=2pt, label=above:$b_1$] at (T) {};
\node[fill=black, circle, inner sep=2pt, label=below:$a_1$] at (B) {};
\node[fill=black, circle, inner sep=2pt, label=left:$a_2$]  at (L) {};
\node[fill=black, circle, inner sep=2pt, label=right:$b_2$] at (R) {};
\node[fill=black, circle, inner sep=2pt, label=-30:$a_3$] at (F) {};
\node[fill=black, circle, inner sep=2pt, label=60:$b_3$] at (K) {};
\node[fill=black, circle, inner sep=2pt, label=210:$c_2$] at (c2) {};
\node[fill=black, circle, inner sep=2pt, label=210:$c_3$] at (c3) {};

\end{scope}

\end{tikzpicture}
}
    \caption{Two graphs whose independence complexes are spheres} 
    \label{fig:Ind(G_3)}
\end{figure}
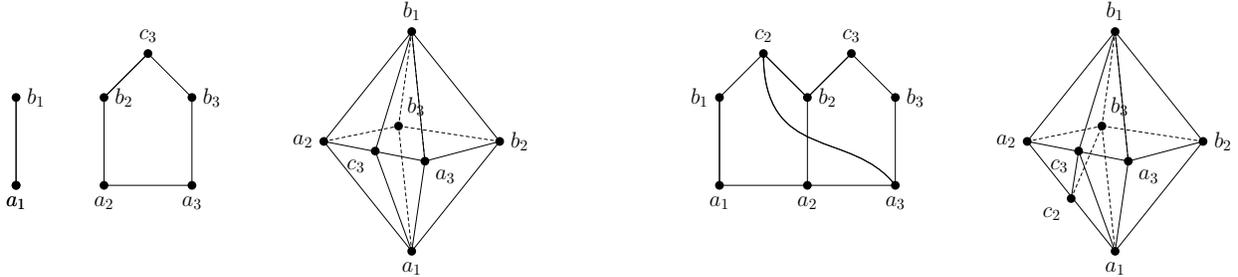

$1$-dimensional simplicial complexes are exactly \emph{simple graphs} (or just ``graphs" from now on). Given a graph $G$, an \emph{independent set} of $G$ is a collection of vertices $S$ of $G$ such that there is no edge in $G$ between elements in $S$. The \emph{independence number}, denoted $\alpha(G)$, is the size of the largest independent set of $G$. The collection of independent sets of $G$ forms a simplicial complex called the \emph{independence complex} of $G$, which we denote by $\Ind(G)$. Simplicial complexes that arise as independence complexes of graphs are \emph{flag complexes}, complexes $\Delta$ where the minimal sets (with respect to inclusion) not in $\Delta$ have size two; in other words, every minimal nonface is a missing edge. \cref{fig:Ind(G_3)} illustrates two graphs and their independence complexes. 

The \emph{induced subcomplex} of $\Delta$ on vertex set $A$ is the complex $\Delta[A]$ on vertex set $A$ with faces $\{\sigma \in \Delta: \sigma \subseteq A\}$. In the case of a graph $G$ with vertex set $V$, we denote the induced subgraph of $G$ on the set $V \setminus A$ as $G - A$.

\begin{remark}[\emph{Well-covered graphs and pure flag complexes}]\label{rmk:well-covered}
    A graph $G$ is \emph{well-covered} if all its maximal independent sets have the same cardinality. Since the maximal independent sets of $G$ correspond to the facets of $\Ind(G)$, this is equivalent to saying $\Ind(G)$ is pure of dimension $\alpha(G)-1$. A well-covered graph is called \emph{$1$-well-covered} if the deletion of any vertex leaves a graph that is also well-covered. In 1979, Staples~\cite{staples1979} introduced the class of $W_n$ graphs, which consists of graphs where every $n$ disjoint independent sets are contained in $n$ disjoint maximum independent sets. Levit and Mandrescu~\cite[Theorem 2.2]{LM2016} then showed that $1$-well-covered graphs without isolated vertices are exactly the graphs in the class $W_2$.

\end{remark}
\subsection{Links, deletions, stars and joins}\label{sec:join}
If $\Delta$ is a simplicial complex and $v$ is a vertex of $\Delta$, we denote by $\lk_{\Delta}(v) = \{\tau \in \Delta \st v \not \in \tau \qand \tau \cup v \in \Delta\}$ the \emph{link} of $v$ in $\Delta$, $\del_\Delta(v) = \{\tau \in \Delta \st v \not \in \tau\}$ the \emph{deletion} of $v$ in $\Delta$, and $\starr_\Delta(v) = \{\tau \in \Delta \st \tau \cup v \in \Delta\}$ the \emph{star} of $v$ in $\Delta$. If there exists a vertex $v$ of $\Delta$ such that $\starr_\Delta(v) = \Delta$, then $\Delta$ is a \emph{cone with apex $v$}, and if $\Delta$ is a cone with apex $v$ for every $v \in \sigma$ for some $\sigma \in \Delta$, we say $\Delta$ is a \emph{cone with $\sigma$}. Note that $\lk_\Delta(v) = \del_\Delta(v) \cap \starr_\Delta(v)$ and $\Delta = \del_\Delta(v) \cup \starr_\Delta(v)$. More generally, the link of a face $\sigma$ of $\Delta$ is the complex $\lk_\Delta(\sigma) = \{\tau \in \Delta \st \tau \cup \sigma \in \Delta \qand \tau \cap \sigma = \emptyset\}$. The coning operation is a special case of the joining operation; given two simplicial complexes $\Delta$ and $\Gamma$ on disjoint vertex sets, the \emph{join} of $\Delta$ and $\Gamma$ is $\Delta \ast \Gamma = \{\tau \sqcup \sigma: \tau \in \Delta, \sigma \in \Gamma\}$. 

For a vertex $v$ of a graph $G$ with vertex set $V$ and edge set $E$, we denote by $N(v) =\{u \in V \st \{u,v\}\in E\}$ the \emph{neighborhood of $v$} and by $N[v] = \{v\} \cup N(v)$ the \emph{closed neighborhood of $v$}. We will often use the following observations:
$$
    \lk_{\Ind(G)}(v) = \Ind(G - N[v]), \quad \del_{\Ind(G)}(v) = \Ind(G - v) \mbox{ and } \starr_{\Ind(G)}(v) = \Ind(G - N(v)).
$$ 
Moreover, if $G$ and $H$ are disjoint graphs, then one can check by the definition of independent set that
\[
\Ind(G \sqcup H) = \Ind(G) \ast \Ind(H). 
\]

\subsection{Cohen--Macaulay complexes, Gorenstein graphs, homology spheres and pseudomanifolds} A simplicial complex $\Delta$ is \emph{Cohen--Macaulay} over a field $\K$ if for every $\sigma \in \Delta$
    $$
        \tilde H_i (\lk_{\Delta}(\sigma); \K) = 0 \qfor i < \dim \lk_\Delta(\sigma), 
    $$
    where $\tilde H_i(\Gamma; \K)$ denotes the reduced $i$-th simplicial homology group of $\Gamma$ over $\K$. If moreover $\Delta$ satisfies $\tilde H_{\dim \lk_\Delta(\sigma)}(\lk_\Delta(\sigma); \K) \cong \K$ for every $\sigma \in \Delta$, then we say $\Delta$ is a \emph{$\K$-homology sphere}. Finally, if $\Delta$ is a homology sphere, or is a cone with $\sigma$ and $\lk_\Delta(\sigma)$ is a homology sphere, then we say $\Delta$ is \emph{Gorenstein}. A graph $G$ for which $\text{Ind}(G)$ is Gorenstein is called a \emph{Gorenstein graph} (for more details see~\cite[II Theorem 5.1]{Stanley}).

A $(d-1)$-dimensional simplicial complex $\Delta$ is a \emph{pseudomanifold} if the following $3$ conditions are satisfied:
    \begin{enumerate}
    \item $\Delta$ is pure;
    \item Given any two facets $F_1, F_2$ of $\Delta$, there exists a sequence of facets $G_0, \dots, G_s$ such that $G_0 = F_1$, $G_s = F_2$ and $|G_i \cap G_{i + 1}| = d-1$. In other words, $\Delta$ is \emph{strongly connected}; and
    \item Every $(d-2)$-dimensional face of $\Delta$ is contained in at most two facets of $\Delta$.\end{enumerate}
    
    The \emph{boundary} of a $(d-1)$-dimensional pseudomanifold $\Delta$ is the collection of $(d-2)$-faces of $\Delta$ contained in exactly one facet. If the boundary of a pseudomanifold is empty, we say it is a pseudomanifold without boundary. 

It is well known that a Cohen--Macaulay complex $\Delta$ is pure and strongly connected, and hence satisfies properties $(1)$ and $(2)$. In particular, a Cohen--Macaulay complex $\Delta$ can only fail to be a pseudomanifold if $\Delta$ has a $(d-2)$-face contained in more than two facets of $\Delta$.

\begin{remark}
    It turns out that the class of $W_2$ graphs (see Remark \ref{rmk:well-covered}) plays a key role in the study of flag homology spheres. It was shown in~\cite[Proposition 3.7]{HoangMinhTrung2016} that a triangle-free graph $G$ without isolated vertices has a Gorenstein independence complex if and only if $G$ is a $W_2$ graph.
\end{remark}

\subsection{Vertex-decomposability, simplicial spheres, and polytopes}
A pure simplicial complex $\Delta$ is \emph{vertex decomposable} if it is a simplex (including $\emptyset$) or there is a vertex $v \in \Delta$ such that the link $\lk_\Delta (v)$ and deletion $\del_\Delta (v)$ are vertex decomposable simplicial complexes.

A \emph{polytope} $P$ is the convex hull of finitely many points in a Euclidean space. Its \emph{dimension} is given by the dimension of the affine hull of $P$. 
For a $d$-polytope $P \subset \mathbb{R}^d$, the proper \emph{faces} of $P$ are the intersections of $P$ with its \emph{supporting hyperplanes}. The polytope $P$ is \emph{simplicial} if all of its $(d-1)$-dimensional faces are simplices. The \emph{boundary complex} of a simplicial $d$-polytope is the simplicial complex consisting of the vertex sets of all faces of dimension less than $d$. We refer the reader to \cite{Ziegler} for more exposition on polytopes. 

A $(d-1)$-dimensional simplicial complex is a \emph{simplicial $(d-1)$-sphere} if its geometric realization in Euclidean space is homeomorphic to a $(d-1)$-sphere. Simplicial spheres are examples of pseudomanifolds and are necessarily homology spheres. The Poincaré homology sphere is a homology sphere but not a simplicial sphere. 

The boundary complex of any simplicial polytope is a simplicial sphere. A simplicial $(d-1)$-sphere is said to be \emph{polytopal} if it can be realized as the boundary complex
of a simplicial $d$-polytope. An example of a nonpolytopal simplicial $3$-sphere is the Barnette sphere \cite{Barnette}.

\subsection{PL spheres and stellar moves.}\label{sec:plandstellar} A simplicial $(d-1)$-sphere is \emph{piecewise linear (PL)} if it is PL homeomorphic to the boundary of a $d$-simplex. Let $\Delta$ be a simplicial complex and $\sigma$ be a face of $\Delta$. Then the \emph{stellar subdivision} of $\Delta$ at $\sigma$ is given by
\begin{enumerate}
    \item first, introducing a new vertex $a$ in the interior of $\sigma$, and
    \item next, replacing $\starr_\Delta (\sigma)$ in $\Delta$ with $a \ast \partial \sigma \ast \lk_\Delta (\sigma)$.
\end{enumerate} The inverse operation is called a \emph{stellar weld}. When $\sigma$ is an edge of $\Delta$, these operations are also called \emph{edge subdivision} (see also \cref{def:edgesubdivision}) and \emph{edge contraction}.

A classical result by Alexander \cite{Alexander1930TheCT} states that two simplicial complexes are PL homeomorphic if and only if they can be connected by a sequence of stellar moves. More recently, Lutz and Nevo \cite{lutznevo16} showed that two flag simplicial spheres are PL homeomorphic if and only if they can be connected by a sequence of edge subdivisions and contractions such that all the complexes in the sequence are flag. The barycentric subdivision of the double suspension of the Poincaré homology sphere is an example of a non-PL flag sphere.

\subsection{\texorpdfstring{$f$-, $h$-, $\gamma$-vectors and $\gamma$-positivity}{f-, h-, gamma-vectors and gamma-positivity}}
\label{sec:fhg}
For a simplicial complex $\Delta$ of dimension $d-1$, we denote by $f_i(\Delta)$ the number of $i$-dimensional faces of $\Delta$. The \emph{$f$-vector of $\Delta$} is the sequence $f(\Delta) = (f_{-1}(\Delta), \dots, f_{d-1}(\Delta))$, and the $f$-polynomial is the generating function of $f(\Delta)$. The \emph{$h$-vector of $\Delta$} is the sequence $h(\Delta) = (h_0(\Delta), \dots, h_{d}(\Delta))$ obtained from $f(\Delta)$ via the polynomial relations 
$$
    \sum_{i = 0}^d h_i(\Delta) t^{d-i} = \sum_{i = 0}^d f_{i-1}(\Delta) (t - 1)^{d - i} .
$$

The \emph{$f$-polynomial} and \emph{$h$-polynomial} of $\Delta$ are respectively $f(\Delta,t):= \sum_{i = 0}^d f_{i-1}(\Delta) t^i$ and $h(\Delta,t):= \sum_{i = 0}^d h_i(\Delta) t^i$. It is well-known that
\[
f(\Delta \ast \Gamma, t) = f(\Delta, t) \cdot f(\Gamma, t) \qand h(\Delta \ast \Gamma, t) = h(\Delta, t) \cdot h(\Gamma, t).
\]
If it is clear from the context, we may omit $\Delta$ in the previous definitions.

When $\Delta$ is a $(d-1)$-dimensional homology sphere, the $h$-vector of $\Delta$ satisfies the Dehn--Sommerville relations: $h_i = h_{d - i}$ for $0 \leq i \leq d$. It is then possible to express $h(\Delta, t)$ as $h(\Delta, t) = \sum_{i = 0}^{\lfloor d/2\rfloor} \gamma_i t^i(1 + t)^{d - 2i}
$. 
The \emph{$\gamma$-vector} of $\Delta$ is the sequence $\gamma(\Delta) = (\gamma_0, \dots, \gamma_{\lfloor d/2 \rfloor})$. If $\gamma_i \geq 0$ for every $i$, we say $h(\Delta, t)$ is \emph{$\gamma$-positive} (or simply $\Delta$ is $\gamma$-positive). In 2005, Gal~\cite{gal2005} conjectured that every flag homology sphere is $\gamma$-positive. As is pointed out in~\cite[Remark 3.1.1]{gal2005} (see also~\cite[Lemma 4.1]{B2004}), if every root of the polynomial $h(\Delta, t)$ is real, then $h(\Delta, t)$ is $\gamma$-positive.

\section{Flag spheres arising from ternary graphs} \label{sec:PropertiesFlagSpheres}

In 1978, Lovász~\cite{lovasz1978} proved Kneser's conjecture by showing connections between topological invariants of certain spaces arising from a graph $G$ and the chromatic number of $G$. Since then, many topological obstructions to graph coloring have been shown. One particular example can be found in the class of \emph{ternary graphs}, that is, graphs without induced cycles of length divisible by three. Two conjectures of Kalai and Meshulam~\cite[Conjectures 1, 2]{blogKalai} would imply there exists a constant $C$, such that the chromatic number of every ternary graph $G$ is bounded above by $C$. The bound on chromatic numbers of ternary graphs was proven by Bonamy, Charbit, and Thomassé in~\cite{BCT2014}, while the following topological characterization was proven by Kim.

\begin{theorem}[{\cite[Theorem 1.1]{kim}}]
\label{t:kim}
    A graph $G$ is ternary if and only if the independence complex of every induced subgraph of $G$ is either contractible or homotopy equivalent to a sphere.
\end{theorem}

In this section, we take a closer look at ternary graphs and their independence complexes.
We prove Theorem \ref{thm:ternarygorenstein} characterizing when the independence complex of a ternary graph is a homology sphere. Then we turn our attention to a family of graphs $\{G_m\}$, which together with the graph $R_3$ give all planar Gorenstein graphs \cite{trung2018characterization}. We show that each $G_m$ is ternary and provide some additional observations about their structure. Finally, we prove that the independence complex of each $G_m$ can be realized as the boundary of a polytope.

For a ternary graph $G$ and a sequence of vertices $v_1, \dots, v_s$ of $G$, a direct consequence of~\cref{t:kim} is that by recursively taking links (or deletions) of the $v_i$, the resulting complex is either contractible or homotopy equivalent to a sphere. This property should be contrasted with the definition of homology spheres. On the one hand, subcomplexes of homology spheres obtained by taking links always have the homology of a sphere of the highest possible dimension. On the other hand, at the homology level, the condition from~\cref{t:kim} includes more subcomplexes than just the ones obtained from links, but it does not require the spheres to be of the highest possible dimension.
In~\cite{ternary} the authors introduce the class of \emph{spherical complexes} as the class of complexes satisfying the topological condition arising from~\cref{t:kim}. The algebraic perspective from~\cite{ternary} led the authors to conjecture a characterization of Gorenstein spherical complexes, which we state for ternary graphs. We give an affirmative answer to~\cite[Question 6.5]{ternary} below.

\begin{lemma}\label{l:topdimCM}
    Let $G$ be a ternary graph and $d = \dim \Ind(G)$. If $\Ind(G)$ is homotopy equivalent to a $d$-dimensional sphere, then $\Ind(G)$ is Cohen--Macaulay.
\end{lemma}

\begin{proof}
Let $\Delta=\Ind(G)$.    By~\cite[Theorem 3.6]{ternary}, we know that for any $\sigma \in \Delta$ either $\lk_\Delta(\sigma)$ has trivial homology or $\lk_\Delta(\sigma)$ has the homology of a $d'$-dimensional sphere, where 
    $$
        d - |\sigma| \leq d' \leq d.
    $$
In particular, $\tilde H_i(\lk_\Delta(\sigma); \K) = 0$ for $i < d - |\sigma|$, and since $\dim \lk_\Delta(\sigma) \leq d - |\sigma|$, we conclude $\Delta$ is Cohen--Macaulay.
\end{proof}

\begin{lemma}\label{l:pseudomanifold}
    Let $G$ be a ternary graph.  If $\Ind(G)$ is Cohen--Macaulay, then $\Ind(G)$ is a pseudomanifold.
\end{lemma}

\begin{proof}
    Since $\Ind(G)$ is Cohen--Macaulay, we conclude it is pure and strongly connected. Let $\sigma$ be a $(d-1)$-dimensional face of $\Ind(G)$. Then $\lk_{\Ind(G)}(\sigma) = \Ind(G - N[\sigma])$, and since $G - N[\sigma]$ is also a ternary graph, we conclude the $0$-dimensional complex $\lk_{\Ind(G)}(\sigma)$ is either contractible or homotopy equivalent to $S^0$. More specifically, $\lk_{\Ind(G)}(\sigma)$ is either one point or two points. In other words, there are at most two vertices of $\Ind(G)$ such that $v \not \in \sigma$ and $\sigma \cup v \in \Ind(G)$. This property, together with the fact that $\Ind(G)$ is pure and strongly connected, implies $\Ind(G)$ is a (Cohen--Macaulay) pseudomanifold.
\end{proof}

\thmA*

\begin{proof}
Let $G$ be a ternary graph without isolated vertices.  Suppose $\Ind(G)$ is a homology sphere. Then $\Ind(G)$ has the same homology as $S^d$ by~\cite[Theorem II.5.1]{Stanley}. \cref{t:kim} then implies $\Ind(G)$ is homotopy equivalent to $S^d$.
    
    Suppose now that $\Ind(G)$ is homotopy equivalent to $S^d$, where $d = \dim \Ind(G)$. By \cref{l:topdimCM} we know $\Ind(G)$ is Cohen--Macaulay, hence by \cref{l:pseudomanifold} we conclude $\Ind(G)$ is a Cohen--Macaulay pseudomanifold.

We now show that $\Ind(G)$ is a Cohen--Macaulay pseudomanifold without boundary. Since $\Ind(G)$ is Cohen--Macaulay and homotopy equivalent to $S^d$, $\tilde H_d (\Ind(G); \Z_2) = \Z_2$. Consider the corresponding chain complex 
\[
\begin{tikzcd}
0 \arrow[r] & C_d \arrow[r, "\partial_d"] & C_{d-1} \arrow[r, "\partial_{d-1}"] & \cdots \arrow[r, "\partial_1"] & C_0 \arrow[r, "\partial_0"] & C_{-1} \arrow[r] & 0.
\end{tikzcd}
\]
Since $\tilde H_d$ is nontrivial, there exists a nonempty set $\{F_i\}$ of facets of $\Ind(G)$ such that $\partial_d(\sum F_i) = 0$. Let $\Delta$ be the subcomplex of $\Ind(G)$ generated by $\{F_i\}$. Since $\partial_d(\sum F_i) = 0$, every ridge of $\Delta$ is contained in an even number of facets, and thus it must be contained in exactly two since $\Ind(G)$ is a pseudomanifold. Thus, $\Delta$ is a $d$-dimensional pseudomanifold without boundary, forcing $\Delta = \Ind(G)$. 

Finally, since $\Ind(G)$ is a Cohen--Macaulay orientable pseudomanifold,~\cite[Theorem II.5.1]{Stanley} directly implies $\Ind(G)$ is Gorenstein.
\end{proof}

\begin{remark}
    Although we state~\cref{thm:ternarygorenstein} for ternary graphs, as it fits the context of this paper, the same proof works for arbitrary spherical complexes (see~\cite[Definition 2.3]{ternary}). 
\end{remark}

In \cite{pinter1995class} a recursive construction of planar $1$-well-covered graphs, $G_m$, $m \in \mathbb{N}$, of girth four is given. 

\begin{definition}
Let $G_1$ be the single edge $a_1b_1$, and $G_2$ be the $5$-cycle $(a_1, b_1, c_2, b_2, a_2, a_1)$. 
For $m \ge 3$, define $G_m$ to be the graph with vertex and edge sets $V_m$ and $E_m$, where
\begin{align*}
V_m &= V(G_m)=V(G_{m-1}) \cup \{a_m, b_m,c_m\} \quad \text{and} \\
E_m &= E(G_m)=E(G_{m-1}) \cup \{a_{m-1}a_m, b_{m-1}c_m, c_{m-1}a_m, a_mb_m, b_mc_m\}.
\end{align*}
\end{definition}

\begin{figure}[H]
    \centering
    \resizebox{.5\textwidth}{!}{
\begin{tikzpicture}[font=\large]
  \useasboundingbox (0,-0.5) rectangle (12,2.7);
\coordinate (ai) at (0, 0);
\coordinate (ai+1) at (2, 0);
\coordinate (ai+2) at (4, 0);
\coordinate (ai+3) at (6, 0);
\coordinate (ai+j-1) at (8, 0);
\coordinate (ai+j) at (10, 0);
\coordinate (ai+j+1) at (12, 0);
\coordinate (bi) at (0, 2);
\coordinate (bi+1) at (2, 2);
\coordinate (bi+2) at (4, 2);
\coordinate (bi+3) at (6, 2);
\coordinate (bi+j-1) at (8, 2);
\coordinate (bi+j) at (10, 2);
\coordinate (bi+j+1) at (12, 2);

\node at (-1,3) (ci) {};
\coordinate (ci+1) at (1, 3);
\coordinate (ci+2) at (3, 3);
\coordinate (ci+3) at (5, 3);
\coordinate (ci+j) at (9, 3);
\coordinate (ci+j+1) at (11, 3);

\draw[thick] (ai) -- (ai+1) -- (ai+2) -- (ai+3);
\draw[thick] (ai+j-1) -- (ai+j) -- (ai+j+1);

\node at (7,1) {\Huge $\cdots$};
\draw[thick] (ai) -- (bi);
\draw[thick] (ai+1) -- (bi+1);
\draw[thick] (ai+2) -- (bi+2);
\draw[thick] (ai+3) -- (bi+3);
\draw[thick] (ai+j-1) -- (bi+j-1);
\draw[thick] (ai+j) -- (bi+j);
\draw[thick] (ai+j+1) -- (bi+j+1);

\draw[thick] (bi) -- (ci+1) -- (bi+1) --(ci+2) -- (bi+2) -- (ci+3) -- (bi+3);
\draw[thick] (bi+j-1)--(ci+j)--(bi+j)--(ci+j+1)--(bi+j+1);

\node[fill=black, circle, inner sep=2pt, label=below:$a_1$] at (ai) {};
\node[fill=black, circle, inner sep=2pt, label=below:$a_{2}$] at (ai+1) {};
\node[fill=black, circle, inner sep=2pt, label=below:$a_3$] at (ai+2) {};
\node[fill=black, circle, inner sep=2pt, label=below:$a_4$] at (ai+3) {};
\node[fill=black, circle, inner sep=2pt, label=below:$a_{m-1}$] at (ai+j-1) {};
\node[fill=black, circle, inner sep=2pt, label=below:$a_m$] at (ai+j) {};
\node[fill=black, circle, inner sep=2pt, label=below:$a_{m+1}$] at (ai+j+1) {};

\node[fill=black, circle, inner sep=2pt, label=-10:$b_1$] at (bi) {};
\node[fill=black, circle, inner sep=2pt, label=-10:$b_2$] at (bi+1) {};
\node[fill=black, circle, inner sep=2pt, label=-10:$b_3$] at (bi+2) {};
\node[fill=black, circle, inner sep=2pt, label=-10:$b_4$] at (bi+3) {};
\node[fill=black, circle, inner sep=2pt, label=-10:$b_{m-1}$] at (bi+j-1) {};
\node[fill=black, circle, inner sep=2pt, label=-10:$b_{m}$] at (bi+j) {};
\node[fill=black, circle, inner sep=2pt, label=-10:$b_{m+1}$] at (bi+j+1) {};

\node[fill=black, circle, inner sep=2pt, label=left:$c_2$] at (ci+1) {};
\node[fill=black, circle, inner sep=2pt, label=left:$c_3$] at (ci+2) {};
\node[fill=black, circle, inner sep=2pt, label=left:$c_4$] at (ci+3) {};
\node[fill=black, circle, inner sep=2pt, label=left:$c_{m}$] at (ci+j) {};
\node[fill=black, circle, inner sep=2pt, label=left:$c_{m+1}$] at (ci+j+1) {};

\draw[thick] (ci+1) .. controls (1,.5) and (3,1.3) .. (ai+2);
\draw[thick] (ci+2) .. controls (3,.5) and (5,1.3) .. (ai+3);
\begin{scope}
    \clip (4,3) rectangle (6.2,0);  
    \draw[thick] (ci+3) .. controls (5,.5) and (7,1.3) .. (ai+j-1);
  \end{scope}
  \begin{scope}
    \clip (7.8,3) rectangle (10,0);  
    \draw[thick] (7,3) .. controls (7,.5) and (9,1.3) .. (ai+j);
  \end{scope}
\draw[thick] (ci+j) .. controls (9,.5) and (11,1.3) .. (ai+j+1);

\end{tikzpicture}}
    \caption{$G_m$}
    \label{fig:G_m_Definition}
\end{figure}

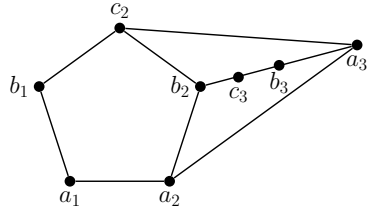
\begin{figure}[H]
  \centering
  \resizebox{.3\linewidth}{!}{%

\begin{tikzpicture}[font=\LARGE]

  \begin{scope}[rotate=-18]
    \node[circle, fill=black, inner sep=0pt, minimum size=9pt, label=180:{$b_1$}] (b1) at (0:-2.5) {};
    \node[circle, fill=black, inner sep=0pt, minimum size=9pt, label=270:{$a_1$}] (a1) at (72:-2.5) {};
    \node[circle, fill=black, inner sep=0pt, minimum size=9pt, label=270:{$a_2$}] (a2) at (144:-2.5) {};
    \node[circle, fill=black, inner sep=0pt, minimum size=9pt, label=180:{$b_2$}] (b2) at (216:-2.5) {};
    \node[circle, fill=black, inner sep=0pt, minimum size=9pt, label=90:{$c_2$}] (c2) at (288:-2.5) {};

    \draw[very thick] (b2)--(c2)--(b1)--(a1)--(a2)--(b2);
  \end{scope}

  \node[circle, fill=black, inner sep=0pt, minimum size=9pt, label=below:{$a_3$}] (a3) at (7,2) {};
  \node[circle, fill=black, inner sep=0pt, minimum size=9pt, label={[label distance=-3pt]below:{$b_3$}}] (b3) at (4.7,1.4) {};
  \node[circle, fill=black, inner sep=0pt, minimum size=9pt, label=below:{$c_3$}] (c3) at (3.5,1.05) {};

  \draw[very thick] (b2) -- (a3);
  \draw[very thick] (a2) -- (a3);
  \draw[very thick] (c2) -- (a3);

\end{tikzpicture}

  }
  \caption{$G_3$ is planar}
  \label{fig:PlanarG_3}
\end{figure}

See Figures \ref{fig:G_m_Definition} and \ref{fig:PlanarG_3} for representations of $G_m$ and $G_3$ respectively. In \cite{hoang2016characterization} this construction is shown to produce Gorenstein graphs, and in \cite{trung2018characterization} it is shown that these are all the triangle-free, connected, planar (see Figure \ref{fig:PlanarG_3}),
Gorenstein graphs without isolated vertices.

\begin{figure}[H]
      \centering
      \resizebox{.2\linewidth}{!}{%
      \begin{tikzpicture}

  \draw[thick] (-2,1) -- (-1,0);
  \draw[thick] (2,1) -- (1,0);
  \draw[thick] (-2,-1) -- (-1,0);
  \draw[thick] (-2,1) -- (-2,-1);
  \draw[thick] (2,1) -- (2,-1);
  \draw[thick] (-2,1) -- (2,1);
  \draw[thick] (-2,-1) -- (2,-1);
  \draw[thick] (2,-1) -- (1,0);
  \draw[thick] (-1,0) -- (1,0);

  \fill (-2,1) circle (2pt);
  \fill (2,1) circle (2pt);
  \fill (-1,0) circle (2pt);
  \fill (1,0) circle (2pt);
  \fill (-2,-1) circle (2pt);
  \fill (2,-1) circle (2pt);

\end{tikzpicture}%
      }
      \caption{The graph $R_3$}
      \label{fig:R3_graph}
    \end{figure}
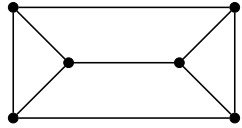

More generally, the following holds.
Here $R_3$ is the graph of a triangular prism (see Figure \ref{fig:R3_graph}).  The independence complex of $R_3$ is the complement of $R_3$, which is a 6-cycle. 
\begin{theorem}[{\cite[Theorem 3.7]{trung2018characterization}}]\label{t:trungcharacterization}
    If a planar graph $G$ is Gorenstein, then $G$ is a disjoint union of graphs of the form $G_m$ and $R_3$. 
    \end{theorem}

We now establish some more properties of $G_m$ for any $m \in \mathbb{N}$, unless otherwise stated. Recall that it was proven in \cite{trung2018characterization} that $G_m$ is triangle-free. In fact, we can show the following stronger statement.
\begin{lemma}\label{l:gmternary}
    $G_m$ is ternary for all $m$.
\end{lemma}
\begin{proof}
    We proceed by induction on $m$.
    Recall that $G_1$ is isomorphic to a single edge $K_2$, and $G_2$ is isomorphic to a $5$-cycle $C_5$. Thus, we immediately see that $G_1$ and $G_2$ are ternary.
    Now assume that $G_{m-1}$ is ternary for some $m \geq 3$.
    Recall that $V(G_m)=V(G_{m-1}) \cup \{a_m, b_m,c_m\}$ and $E(G_m)=E(G_{m-1}) \cup \{a_{m-1}a_m, b_{m-1}c_m, c_{m-1}a_m, a_mb_m, b_mc_m\}$.
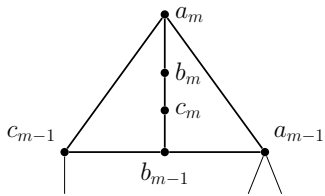
    \begin{figure}[H]
        \centering
        \resizebox{0.27\textwidth}{!}{\begin{tikzpicture}[every node/.style={font=\Large}]

\node[fill=black, circle, inner sep=2pt, label=130:$c_{m-1}$] (a) at (-.4, 0) {};
\node[fill=black, circle, inner sep=2pt, label=below:$b_{m-1}$] (v) at (2, 0) {};
\node[fill=black, circle, inner sep=2pt, label=50:$a_{m-1}$] (b) at (4.4, 0) {};

\node[fill=black, circle, inner sep=2pt, label=right:$a_m$] (y) at (2, 3.3) {};
\node[fill=black, circle, inner sep=2pt, label=right:$b_m$] (x) at (2, 1.9) {};
\node[fill=black, circle, inner sep=2pt, label=right:$c_m$] (w) at (2, 1) {};

\draw[black, very thick] (a) -- (y);
\draw[black, very thick] (b) -- (y);
\draw[very thick] (a) -- (v);
\draw[very thick] (b) -- (v);
\draw[black, very thick] (y) -- (x) -- (w) -- (v);

\draw (a) -- ++(0, -1);
\draw (b) -- ++(-.4, -1);
\draw (b) -- ++(.4, -1);

\end{tikzpicture}}
        \label{fig:GmInduction}
    \end{figure}
    \noindent Note that all induced cycles in $G_{m-1}$ are induced cycles in $G_m$. As $G_{m-1}$ is ternary, none of these induced cycles have length $3k$ for any $k \in \mathbb{N}$.
    Additionally, $G_m$ has induced $5$-cycles $(a_{m-1}, b_{m-1}, c_m, b_m, a_m, a_{m-1})$ and $(b_{m-1}, c_{m-1}, a_m, b_m, c_m, b_{m-1})$, and an induced $4$-cycle $(a_{m-1}, b_{m-1}, c_{m-1}, a_m, a_{m-1})$.
    Finally, we can build induced cycles in $G_m$ from induced cycles in $G_{m-1}$ as follows. Observe that as $\deg_{G_{m-1}}(b_{m-1})=2$, a cycle in $G_{m-1}$ contains the edge $a_{m-1}b_{m-1}$ if and only if it contains the edge $b_{m-1}c_{m-1}$. Let $C$ be such a cycle in $G_{m-1}$, in particular, $\len(C) \neq 3k$ for any $k \in \mathbb{N}$. Replacing the path $(a_{m-1},b_{m-1},c_{m-1})$ in $C$ by the path $(a_{m-1},a_m,c_{m-1})$ results in an induced cycle $C'$ in $G_m$ of length $\len(C') = \len(C) \neq 3k$. 
    We can now conclude that none of the induced cycles in $G_m$ have length $3k$ for any $k \in \mathbb{N}$, and thus $G_m$ is ternary as desired.
\end{proof}

\begin{remark}\label{rem:R3}
    In view of~\cref{t:trungcharacterization,l:gmternary}, we will often state our results for planar ternary graphs. For results we state with this assumption, there is a less clean version including the graph $R_3$, which is the only nonternary connected planar Gorenstein graph. For the sake of brevity, we omit these modified versions.
\end{remark}

Having excluded the existence of induced cycles of length divisible by three, we now investigate the existence and interactions of other cycles in $G_m$.

\begin{proposition}
    \label{prop:cycles_gm} Let $m \ge 1$. Then $G_m$ has the following properties:
    \begin{enumerate}[(a)]
        \item Any two distinct 4-cycles in $G_m$ are edge-disjoint. 
        \item $G_m$ has no $6$-cycles.
        \item Any pair consisting of a $4$-cycle and a $5$-cycle in $G_m$ is either edge-disjoint or intersects in a $2$-path. 
    \end{enumerate}
\end{proposition}

\begin{proof}
    \text{}
    \begin{enumerate}[(a)]
        \item By construction of $G_m$, all of its $4$-cycles are of the form $(a_{t-1}, a_{t}, c_{t-1}, b_{t-1}, a_{t-1})$ for some $t \geq 3$. Thus, any two distinct $4$-cycles are edge-disjoint.
        \item Suppose $G_m$ had a $6$-cycle, say $(x_1,x_2,\ldots,x_6,x_1)$. As $G_m$ is ternary, there exists an edge $x_ix_{i+3}$ (indices mod 6)
        for some $i \in [6]$. However, this gives two distinct $4$-cycles sharing one edge, which do not exist by part (a). Hence, $G_m$ has no $6$-cycles.
        \item Recall that every $4$-cycle in $G_m$ is of the form $(a_{t-1},a_{t},c_{t-1},b_{t-1},a_{t-1})$ for some $t \geq 3$. By construction of $G_m$, every $5$-cycle is of the form $(a_{t-1},a_t,b_t,c_t,b_{t-1},a_{t-1})$ for some $t \geq 2$, $(a_{t-2},a_{t-1},a_t,c_{t-1},b_{t-2},a_{t-2})$ for some $t \geq 3$, $(a_t, b_t, c_t, b_{t-1}, c_{t-1}, a_t)$ for some $t \geq 3$, or $(a_{t-1},a_t,c_{t-1},b_{t-2},c_{t-2},a_{t-1})$ for some $t \geq 4$. If a $4$-cycle and a $5$-cycle are not edge-disjoint, it is now straightforward to see that the intersection is a $2$-path.
    \end{enumerate}
\end{proof}

We now turn our attention to induced paths in $G_m$. We find that unless two vertices are close to each other, there exist multiple induced paths of different lengths between them. 
\begin{lemma}\label{l:paths_mod3}
    Let $x, y \in V(G_m)$. If $\dist(x,y) \geq 3$, then there exist induced paths of length $0, 1,$ and $2 \mod 3$ between $x$ and $y$ in $G_m$.
\end{lemma}

\begin{proof}
\begin{figure}[htp]
  \centering
  \captionsetup[subfigure]{labelformat=empty, labelsep=none}
  \begin{subfigure}[t]{0.48\linewidth}
    \centering
    \resizebox{\linewidth}{!}{

\begin{tikzpicture}[
  every node/.style={font=\Large},
  every label/.append style={font=\Large},
  baseline=(current bounding box.north) 
]
  \path[use as bounding box] (-1.4,-1.0) rectangle (16.4,4.3);

  \node at (7.5,1)       {\Huge $\cdots$};

  \coordinate (ak)       at ( 0,0);
  \coordinate (ak+1)     at ( 2,0);
  \coordinate (ak+2)     at ( 4,0);
  \coordinate (ak+3)     at ( 6,0);
  \coordinate (ak+j-1)   at ( 9,0);
  \coordinate (ak+j)     at (11,0);
  \coordinate (ak+j+1)   at (13,0);

  \coordinate (bk)       at ( 0,2);
  \coordinate (bk+1)     at ( 2,2);
  \coordinate (bk+2)     at ( 4,2);
  \coordinate (bk+3)     at ( 6,2);
  \coordinate (bk+j-1)   at ( 9,2);
  \coordinate (bk+j)     at (11,2);
  \coordinate (bk+j+1)   at (13,2);

  \coordinate (ck)       at ( -1,3);
  \coordinate (ck+1)     at ( 1,3);
  \coordinate (ck+2)     at ( 3,3);
  \coordinate (ck+3)     at ( 5,3);
  \coordinate (ck+j)     at (10,3);
  \coordinate (ck+j+1)   at (12,3);


  \draw[gray!60,line width=3mm]
    (ak)  -- (ak+3)
    (ak+j-1) -- (ak+j);

  \draw[thick]
    (ak+j)  -- (ak+j+1);

  \draw[thick]
    (ak)      -- (bk)
    (ak+1)    -- (bk+1)
    (ak+2)    -- (bk+2)
    (ak+3)    -- (bk+3)
    (ak+j-1)  -- (bk+j-1)
    (ak+j)    -- (bk+j)
    (ak+j+1)  -- (bk+j+1);

  \draw[thick]
    (bk)      -- (ck+1) -- (bk+1) -- (ck+2) -- (bk+2) -- (ck+3) -- (bk+3)
    (bk+j-1)  -- (ck+j)  -- (bk+j)  -- (ck+j+1) -- (bk+j+1);

  \begin{scope}
    \clip (-.3,0) rectangle (2.1,3.1);
    \draw[thick]
      (bk) -- (ck)
      (ck) .. controls (-1,.5) and (1,1.3) .. (ak+1)
      (ak) -- (-1,0);
  \end{scope}

  \draw[thick]
    (ck+1) .. controls (1,.5)  and (3,1.3)  .. (ak+2)
    (ck+2) .. controls (3,.5)  and (5,1.3)  .. (ak+3);

  \begin{scope}
    \clip (4,3) rectangle (6.3,-.1);
    \draw[thick]
      (ck+3) .. controls (5,.5) and (7,1.3) .. (ak+j-1);
    \draw[thick] (ak+3) -- (7,0);
    \draw[gray!60,line width=.8mm] (ak+3) -- (7,0);
    \draw[thick] (bk+3) -- (7,3);
  \end{scope}

  \begin{scope}
    \clip (8.7,3) rectangle (11,-.1);
    \draw[thick]
      (7,3) .. controls (8,.5)   and (10,1.3) .. (ak+j)
      (bk+j-1) -- (8,3);
    \draw[thick] (ak+j-1) -- (8,0);
    \draw[gray!60,line width=.8mm] (ak+j-1) -- (8,0);
  \end{scope}

  \draw[thick]
    (ck+j) .. controls (10,.5) and (12,1.3) .. (ak+j+1);

  \begin{scope}
    \clip (12,3) rectangle (13.3,-.1);
    \draw[thick]
      (ck+j+1) .. controls (12,.5) and (14,1.3) .. (16,0)
      (bk+j+1) -- (14,3)
      (ak+j+1) -- (14,0);
  \end{scope}

  \node[fill=black,circle,inner sep=2pt,label=below:$a_k$]         at (ak)       {};
  \node[fill=black,circle,inner sep=2pt,label=below:$a_{k+1}$]     at (ak+1)     {};
  \node[fill=black,circle,inner sep=2pt,label=below:$a_{k+2}$]     at (ak+2)     {};
  \node[fill=black,circle,inner sep=2pt,label=below:$a_{k+3}$]     at (ak+3)     {};
  \node[fill=black,circle,inner sep=2pt,label=below:$a_{k+l-1}$]   at (ak+j-1)   {};
  \node[fill=black,circle,inner sep=2pt,label=below:$a_{k+l}$]     at (ak+j)     {};
  \node[fill=black,circle,inner sep=2pt,label=below:$a_{k+l+1}$]   at (ak+j+1)   {};

  \node[fill=black,circle,inner sep=2pt,label=-10:$b_k$]          at (bk)       {};
  \node[fill=black,circle,inner sep=2pt,label=-10:$b_{k+1}$]      at (bk+1)     {};
  \node[fill=black,circle,inner sep=2pt,label=-10:$b_{k+2}$]      at (bk+2)     {};
  \node[fill=black,circle,inner sep=2pt,label=-10:$b_{k+3}$]      at (bk+3)     {};
  \node[fill=black,circle,inner sep=2pt,label=-10:$b_{k+l-1}$]    at (bk+j-1)   {};
  \node[fill=black,circle,inner sep=2pt,label=-10:$b_{k+l}$]      at (bk+j)     {};
  \node[fill=black,circle,inner sep=2pt,label=-10:$b_{k+l+1}$]    at (bk+j+1)   {};

  \node[fill=black,circle,inner sep=2pt,label=above:$c_{k+1}$]    at (ck+1)     {};
  \node[fill=black,circle,inner sep=2pt,label=above:$c_{k+2}$]    at (ck+2)     {};
  \node[fill=black,circle,inner sep=2pt,label=above:$c_{k+3}$]    at (ck+3)     {};
  \node[fill=black,circle,inner sep=2pt,label=above:$c_{k+l}$]    at (ck+j)     {};
  \node[fill=black,circle,inner sep=2pt,label=above:$c_{k+l+1}$]  at (ck+j+1)   {};

\end{tikzpicture}
}
    \caption{Type (1,a)}\label{fig:1a}
  \end{subfigure}\hfill
  \begin{subfigure}[t]{0.48\linewidth}
    \centering
    \resizebox{\linewidth}{!}{\begin{tikzpicture}[
  every node/.style={font=\Large},
  every label/.append style={font=\Large},
  baseline=(current bounding box.north) 
]
  \path[use as bounding box] (-1.4,-1.0) rectangle (16.4,4.3);

  \node at (7.5,1)  {\Huge $\cdots$};

  \coordinate (ai)       at ( 0,0);
  \coordinate (ai+1)     at ( 2,0);
  \coordinate (ai+2)     at ( 4,0);
  \coordinate (ai+3)     at ( 6,0);
  \coordinate (ai+j-1)   at ( 9,0);
  \coordinate (ai+j)     at (11,0);
  \coordinate (ai+j+1)   at (13,0);

  \coordinate (bi)       at ( 0,2);
  \coordinate (bi+1)     at ( 2,2);
  \coordinate (bi+2)     at ( 4,2);
  \coordinate (bi+3)     at ( 6,2);
  \coordinate (bi+j-1)   at ( 9,2);
  \coordinate (bi+j)     at (11,2);
  \coordinate (bi+j+1)   at (13,2);

  \node       at (-1,3) (ci)     {};
  \coordinate (ci+1)   at ( 1,3);
  \coordinate (ci+2)   at ( 3,3);
  \coordinate (ci+3)   at ( 5,3);
  \coordinate (ci+j)   at ( 10,3);
  \coordinate (ci+j+1) at (12,3);

  \draw[gray!60,line width=3mm]
    (bi) -- (ci+1) -- (bi+1) -- (ci+2) -- (bi+2) -- (ci+3) -- (bi+3);

  \draw[gray!60,line width=3mm]
     (bi+j-1) -- (ci+j) --  (bi+j);

  \draw[thick] (ai)      -- (ai+1) -- (ai+2) -- (ai+3);
  \draw[thick] (ai+j-1)  -- (ai+j)  -- (ai+j+1);

  \draw[thick] (ai)      -- (bi);
  \draw[thick] (ai+1)    -- (bi+1);
  \draw[thick] (ai+2)    -- (bi+2);
  \draw[thick] (ai+3)    -- (bi+3);
  \draw[thick] (ai+j-1)  -- (bi+j-1);
  \draw[thick] (ai+j)    -- (bi+j);
  \draw[thick] (ai+j+1)  -- (bi+j+1);

  \draw[thick] (bi+j)  -- (ci+j+1) -- (bi+j+1);

  \begin{scope}
    \clip (-.3,0) rectangle (2.1,3.1);
    \draw[thick] (bi) -- (ci)
                  (ci) .. controls (-1,.5) and (1,1.3) .. (ai+1);
    \draw[thick] (ai)--(-1,0);
  \end{scope}

  \draw[thick]
    (ci+1) .. controls (1,.5) and (3,1.3) .. (ai+2)
    (ci+2) .. controls (3,.5) and (5,1.3) .. (ai+3);

  \begin{scope}
    \clip (4,3) rectangle (6.3,-.1);
    \draw[thick] (ci+3) .. controls (5,.5) and (7,1.3) .. (ai+j-1);
    \draw[gray!60,line width=3mm] (bi+3)    -- (7,3);
    \draw[thick] (ai+3)    -- (7,0);

  \end{scope}

  \begin{scope}
    \clip (8.7,3) rectangle (11,-.1);
    \draw[thick] (7,3) .. controls (8,.5) and (10,1.3) .. (ai+j);
    \draw[gray!60,line width=3mm] (bi+j-1) -- (8,3);
    \draw[thick] (ai+j-1) -- (8,0);

  \end{scope}

  \draw[thick]
    (ci+j) .. controls (10,.5) and (12,1.3) .. (ai+j+1);

  \begin{scope}
    \clip (12,3) rectangle (13.3,-.1);
    \draw[thick] (ci+j+1) .. controls (12,.5) and (14,1.3) .. (16,0);
    \draw[thick] (bi+j+1)--(14,3);
    \draw[thick] (ai+j+1)--(14,0);
  \end{scope}

  \node[fill=black,circle,inner sep=2pt,label=below:$a_k$]           at (ai)       {};
  \node[fill=black,circle,inner sep=2pt,label=below:$a_{k+1}$]       at (ai+1)     {};
  \node[fill=black,circle,inner sep=2pt,label=below:$a_{k+2}$]       at (ai+2)     {};
  \node[fill=black,circle,inner sep=2pt,label=below:$a_{k+3}$]       at (ai+3)     {};
  \node[fill=black,circle,inner sep=2pt,label=below:$a_{k+l-1}$]     at (ai+j-1)   {};
  \node[fill=black,circle,inner sep=2pt,label=below:$a_{k+l}$]       at (ai+j)     {};
  \node[fill=black,circle,inner sep=2pt,label=below:$a_{k+l+1}$]     at (ai+j+1)   {};

  \node[fill=black,circle,inner sep=2pt,label=-10:$b_k$]            at (bi)       {};
  \node[fill=black,circle,inner sep=2pt,label=-10:$b_{k+1}$]        at (bi+1)     {};
  \node[fill=black,circle,inner sep=2pt,label=-10:$b_{k+2}$]        at (bi+2)     {};
  \node[fill=black,circle,inner sep=2pt,label=-10:$b_{k+3}$]        at (bi+3)     {};
  \node[fill=black,circle,inner sep=2pt,label=-10:$b_{k+l-1}$]      at (bi+j-1)   {};
  \node[fill=black,circle,inner sep=2pt,label=-10:$b_{k+l}$]        at (bi+j)     {};
  \node[fill=black,circle,inner sep=2pt,label=-10:$b_{k+l+1}$]      at (bi+j+1)   {};

  \node[fill=black,circle,inner sep=2pt,label=above:$c_{k+1}$]       at (ci+1)     {};
  \node[fill=black,circle,inner sep=2pt,label=above:$c_{k+2}$]       at (ci+2)     {};
  \node[fill=black,circle,inner sep=2pt,label=above:$c_{k+3}$]      at (ci+3)     {};
  \node[fill=black,circle,inner sep=2pt,label=above:$c_{k+l}$]      at (ci+j)     {};
  \node[fill=black,circle,inner sep=2pt,label=above:$c_{k+l+1}$]    at (ci+j+1)   {};

\end{tikzpicture}
}
    \caption{Type (1,b)}\label{fig:1b}
  \end{subfigure}

  \vspace{1em}

  \begin{subfigure}[t]{0.48\linewidth}
    \centering
    \resizebox{\linewidth}{!}{\begin{tikzpicture}[
  every node/.style={font=\Large},
  every label/.append style={font=\Large},
  baseline=(current bounding box.north) 
]
  \path[use as bounding box] (-1.4,-1.0) rectangle (16.4,4.3);

  \node at (7.5,1)       {\Huge $\cdots$};

  \coordinate (ak)       at ( 0,0);
  \coordinate (ak+1)     at ( 2,0);
  \coordinate (ak+2)     at ( 4,0);
  \coordinate (ak+3)     at ( 6,0);
  \coordinate (ak+j-1)   at ( 9,0);
  \coordinate (ak+j)     at (11,0);
  \coordinate (ak+j+1)   at (13,0);

  \coordinate (bk)       at ( 0,2);
  \coordinate (bk+1)     at ( 2,2);
  \coordinate (bk+2)     at ( 4,2);
  \coordinate (bk+3)     at ( 6,2);
  \coordinate (bk+j-1)   at ( 9,2);
  \coordinate (bk+j)     at (11,2);
  \coordinate (bk+j+1)   at (13,2);

  \coordinate (ck)       at ( -1,3);
  \coordinate (ck+1)     at ( 1,3);
  \coordinate (ck+2)     at ( 3,3);
  \coordinate (ck+3)     at ( 5,3);
  \coordinate (ck+j)     at (10,3);
  \coordinate (ck+j+1)   at (12,3);

  \draw[gray!60, line width = 3mm]
    (ak) -- (bk) -- (ck+1)
    (ck+1) .. controls (1,.5) and (3,1.3) .. (ak+2) -- (ak+3)
    (ak+j-1) -- (ak+j);

  \draw[thick]
    (ak)      -- (ak+2)
    (ak+j)  -- (ak+j+1);

  \draw[thick]
    (ak+1)    -- (bk+1)
    (ak+2)    -- (bk+2)
    (ak+3)    -- (bk+3)
    (ak+j-1)  -- (bk+j-1)
    (ak+j)    -- (bk+j)
    (ak+j+1)  -- (bk+j+1);

  \draw[thick]
    (ck+1) -- (bk+1) -- (ck+2) -- (bk+2) -- (ck+3) -- (bk+3)
    (bk+j-1)  -- (ck+j)  -- (bk+j)  -- (ck+j+1) -- (bk+j+1);

  \begin{scope}
    \clip (-.3,0) rectangle (2.1,3.1);
    \draw[thick]
      (bk) -- (ck)
      (ck) .. controls (-1,.5) and (1,1.3) .. (ak+1)
      (ak) -- (-1,0);
  \end{scope}

  \draw[thick]
    (ck+2) .. controls (3,.5)  and (5,1.3)  .. (ak+3);

  \begin{scope}
    \clip (4,3) rectangle (6.3,-.1);
    \draw[thick]
      (ck+3) .. controls (5,.5) and (7,1.3) .. (ak+j-1);
      \draw[thick](ak+3) -- (7,0);
    \draw[gray!60, line width = 3mm] (ak+3) -- (7,0);
    \draw[thick] (bk+3) -- (7,3);
  \end{scope}

  \begin{scope}
    \clip (8.7,3) rectangle (11,-.1);
    \draw[thick]
      (7,3) .. controls (8,.5)   and (10,1.3) .. (ak+j)
      (bk+j-1) -- (8,3);
    \draw[thick] (ak+j-1) -- (8,0);
    \draw[gray!60, line width = 3mm] (ak+j-1) -- (8,0);

  \end{scope}

  \draw[thick]
    (ck+j) .. controls (10,.5) and (12,1.3) .. (ak+j+1);

  \begin{scope}
    \clip (12,3) rectangle (13.3,-.1);
    \draw[thick]
      (ck+j+1) .. controls (12,.5) and (14,1.3) .. (16,0)
      (bk+j+1) -- (14,3)
      (ak+j+1) -- (14,0);
  \end{scope}

  \node[fill=black,circle,inner sep=2pt,label=below:$a_k$]         at (ak)       {};
  \node[fill=black,circle,inner sep=2pt,label=below:$a_{k+1}$]     at (ak+1)     {};
  \node[fill=black,circle,inner sep=2pt,label=below:$a_{k+2}$]     at (ak+2)     {};
  \node[fill=black,circle,inner sep=2pt,label=below:$a_{k+3}$]     at (ak+3)     {};
  \node[fill=black,circle,inner sep=2pt,label=below:$a_{k+l-1}$]   at (ak+j-1)   {};
  \node[fill=black,circle,inner sep=2pt,label=below:$a_{k+l}$]     at (ak+j)     {};
  \node[fill=black,circle,inner sep=2pt,label=below:$a_{k+l+1}$]   at (ak+j+1)   {};

  \node[fill=black,circle,inner sep=2pt,label=-10:$b_k$]          at (bk)       {};
  \node[fill=black,circle,inner sep=2pt,label=-10:$b_{k+1}$]      at (bk+1)     {};
  \node[fill=black,circle,inner sep=2pt,label=-10:$b_{k+2}$]      at (bk+2)     {};
  \node[fill=black,circle,inner sep=2pt,label=-10:$b_{k+3}$]      at (bk+3)     {};
  \node[fill=black,circle,inner sep=2pt,label=-10:$b_{k+l-1}$]    at (bk+j-1)   {};
  \node[fill=black,circle,inner sep=2pt,label=-10:$b_{k+l}$]      at (bk+j)     {};
  \node[fill=black,circle,inner sep=2pt,label=-10:$b_{k+l+1}$]    at (bk+j+1)   {};

  \node[fill=black,circle,inner sep=2pt,label=above:$c_{k+1}$]    at (ck+1)     {};
  \node[fill=black,circle,inner sep=2pt,label=above:$c_{k+2}$]    at (ck+2)     {};
  \node[fill=black,circle,inner sep=2pt,label=above:$c_{k+3}$]    at (ck+3)     {};
  \node[fill=black,circle,inner sep=2pt,label=above:$c_{k+l}$]    at (ck+j)     {};
  \node[fill=black,circle,inner sep=2pt,label=above:$c_{k+l+1}$]  at (ck+j+1)   {};

\end{tikzpicture}
}
    \caption{Type (2,a)}\label{fig:2a}
  \end{subfigure}\hfill
  \begin{subfigure}[t]{0.48\linewidth}
    \centering
    \resizebox{\linewidth}{!}{\begin{tikzpicture}[
  every node/.style={font=\Large},
  every label/.append style={font=\Large},
  baseline=(current bounding box.north) 
]
  \path[use as bounding box] (-1.4,-1.0) rectangle (16.4,4.3);

  \node at (7.5,1)       {\Huge $\cdots$};

  \coordinate (ak)       at ( 0,0);
  \coordinate (ak+1)     at ( 2,0);
  \coordinate (ak+2)     at ( 4,0);
  \coordinate (ak+3)     at ( 6,0);
  \coordinate (ak+j-1)   at ( 9,0);
  \coordinate (ak+j)     at (11,0);
  \coordinate (ak+j+1)   at (13,0);

  \coordinate (bk)       at ( 0,2);
  \coordinate (bk+1)     at ( 2,2);
  \coordinate (bk+2)     at ( 4,2);
  \coordinate (bk+3)     at ( 6,2);
  \coordinate (bk+j-1)   at ( 9,2);
  \coordinate (bk+j)     at (11,2);
  \coordinate (bk+j+1)   at (13,2);

  \coordinate (ck)       at ( -1,3);
  \coordinate (ck+1)     at ( 1,3);
  \coordinate (ck+2)     at ( 3,3);
  \coordinate (ck+3)     at ( 5,3);
  \coordinate (ck+j)     at (10,3);
  \coordinate (ck+j+1)   at (12,3);

  \draw[gray!60,line width=3mm]
    (bk) -- (ck+1)
    (ck+1) .. controls (1,.5) and (3,1.3) .. (ak+2)
    (ai+2) -- (bi+2) -- (ci+3) -- (bi+3)
    (bk+j-1) -- (ck+j);

  \draw[thick]
    (ak)      -- (ak+1) -- (ak+2) -- (ak+3)
    (ak+j-1)  -- (ak+j)  -- (ak+j+1);

  \draw[thick]
    (ak)      -- (bk)
    (ak+1)    -- (bk+1)
    (ak+3)    -- (bk+3)
    (ak+j-1)  -- (bk+j-1)
    (ak+j)    -- (bk+j)
    (ak+j+1)  -- (bk+j+1);

  \draw[thick]
    (ck+1) -- (bk+1) -- (ck+2) -- (bk+2);
    \draw[thick] (ck+j)  -- (bk+j)  -- (ck+j+1) -- (bk+j+1);

  \begin{scope}
    \clip (-.3,0) rectangle (2.1,3.1);
    \draw[thick]
      (bk) -- (ck)
      (ck) .. controls (-1,.5) and (1,1.3) .. (ak+1)
      (ak) -- (-1,0);
  \end{scope}

  \draw[thick]
    (ck+2) .. controls (3,.5)  and (5,1.3)  .. (ak+3);

  \begin{scope}
    \clip (4,3) rectangle (6.3,-.1);
    \draw[thick]
      (ck+3) .. controls (5,.5) and (7,1.3) .. (ak+j-1);
      \draw[thick](ak+3) -- (7,0);
    \draw[gray!60,line width=3mm] (bk+3) -- (7,3);
  \end{scope}

  \begin{scope}
    \clip (8.7,3) rectangle (11,-.1);
    \draw[thick]
      (7,3) .. controls (8,.5)   and (10,1.3) .. (ak+j)
      (bk+j-1) -- (8,3);
    \draw[thick] (ak+j-1) -- (8,0);
    \draw[gray!60,line width=3mm] (bk+j-1) -- (8,3);

  \end{scope}

  \draw[thick]
    (ck+j) .. controls (10,.5) and (12,1.3) .. (ak+j+1);

  \begin{scope}
    \clip (12,3) rectangle (13.3,-.1);
    \draw[thick]
      (ck+j+1) .. controls (12,.5) and (14,1.3) .. (16,0)
      (bk+j+1) -- (14,3)
      (ak+j+1) -- (14,0);
  \end{scope}

  \node[fill=black,circle,inner sep=2pt,label=below:$a_k$]         at (ak)       {};
  \node[fill=black,circle,inner sep=2pt,label=below:$a_{k+1}$]     at (ak+1)     {};
  \node[fill=black,circle,inner sep=2pt,label=below:$a_{k+2}$]     at (ak+2)     {};
  \node[fill=black,circle,inner sep=2pt,label=below:$a_{k+3}$]     at (ak+3)     {};
  \node[fill=black,circle,inner sep=2pt,label=below:$a_{k+l-1}$]   at (ak+j-1)   {};
  \node[fill=black,circle,inner sep=2pt,label=below:$a_{k+l}$]     at (ak+j)     {};
  \node[fill=black,circle,inner sep=2pt,label=below:$a_{k+l+1}$]   at (ak+j+1)   {};

  \node[fill=black,circle,inner sep=2pt,label=-10:$b_k$]          at (bk)       {};
  \node[fill=black,circle,inner sep=2pt,label=-10:$b_{k+1}$]      at (bk+1)     {};
  \node[fill=black,circle,inner sep=2pt,label=-10:$b_{k+2}$]      at (bk+2)     {};
  \node[fill=black,circle,inner sep=2pt,label=-10:$b_{k+3}$]      at (bk+3)     {};
  \node[fill=black,circle,inner sep=2pt,label=-10:$b_{k+l-1}$]    at (bk+j-1)   {};
  \node[fill=black,circle,inner sep=2pt,label=-10:$b_{k+l}$]      at (bk+j)     {};
  \node[fill=black,circle,inner sep=2pt,label=-10:$b_{k+l+1}$]    at (bk+j+1)   {};

  \node[fill=black,circle,inner sep=2pt,label=above:$c_{k+1}$]    at (ck+1)     {};
  \node[fill=black,circle,inner sep=2pt,label=above:$c_{k+2}$]    at (ck+2)     {};
  \node[fill=black,circle,inner sep=2pt,label=above:$c_{k+3}$]    at (ck+3)     {};
  \node[fill=black,circle,inner sep=2pt,label=above:$c_{k+l}$]    at (ck+j)     {};
  \node[fill=black,circle,inner sep=2pt,label=above:$c_{k+l+1}$]  at (ck+j+1)   {};

\end{tikzpicture}
}
    \caption{Type (2,b)}\label{fig:2b}
  \end{subfigure}

  \vspace{1em}

  \begin{subfigure}[t]{0.48\linewidth}
    \centering
    \resizebox{\linewidth}{!}{\begin{tikzpicture}[
  every node/.style={font=\Large},
  every label/.append style={font=\Large},
  baseline=(current bounding box.north) 
]
  \path[use as bounding box] (-1.4,-1.0) rectangle (16.4,4.3);

  \node at (7.5,1)  {\Huge $\cdots$};

  \coordinate (ai)       at ( 0,0);
  \coordinate (ai+1)     at ( 2,0);
  \coordinate (ai+2)     at ( 4,0);
  \coordinate (ai+3)     at ( 6,0);
  \coordinate (ai+j-1)   at ( 9,0);
  \coordinate (ai+j)     at (11,0);
  \coordinate (ai+j+1)   at (13,0);

  \coordinate (bi)       at ( 0,2);
  \coordinate (bi+1)     at ( 2,2);
  \coordinate (bi+2)     at ( 4,2);
  \coordinate (bi+3)     at ( 6,2);
  \coordinate (bi+j-1)   at ( 9,2);
  \coordinate (bi+j)     at (11,2);
  \coordinate (bi+j+1)   at (13,2);

  \node       at (-1,3) (ci)     {};
  \coordinate (ci+1)   at ( 1,3);
  \coordinate (ci+2)   at ( 3,3);
  \coordinate (ci+3)   at ( 5,3);
  \coordinate (ci+j)   at ( 10,3);
  \coordinate (ci+j+1) at (12,3);

  \draw[gray!60,line width=3mm]
    (ai) -- (bi) -- (ci+1) -- (bi+1) -- (ci+2);

  \draw[gray!60,line width=3mm]
    (ci+2) .. controls (3,.5) and (5,1.3) .. (ai+3);

  \draw[gray!60,line width=3mm]
    (ai+j-1) -- (ai+j);

  \draw[thick] (ai)      -- (ai+1) -- (ai+2) -- (ai+3);
  \draw[thick]  (ai+j)  -- (ai+j+1);

  \draw[thick] (ai+1)    -- (bi+1);
  \draw[thick] (ai+2)    -- (bi+2);
  \draw[thick] (ai+3)    -- (bi+3);
  \draw[thick] (ai+j-1)  -- (bi+j-1);
  \draw[thick] (ai+j)    -- (bi+j);
  \draw[thick] (ai+j+1)  -- (bi+j+1);

  \draw[thick] (ci+2) -- (bi+2) -- (ci+3) -- (bi+3);
  \draw[thick] (bi+j-1)  -- (ci+j)  -- (bi+j)  -- (ci+j+1) -- (bi+j+1);

  \begin{scope}
    \clip (-.3,0) rectangle (2.1,3.1);
    \draw[thick] (bi) -- (ci)
                  (ci) .. controls (-1,.5) and (1,1.3) .. (ai+1);
    \draw[thick] (ai)--(-1,0);
  \end{scope}

  \draw[thick]
    (ci+1) .. controls (1,.5) and (3,1.3) .. (ai+2);

  \begin{scope}
    \clip (4,3) rectangle (6.3,-.1);
    \draw[thick] (ci+3) .. controls (5,.5) and (7,1.3) .. (ai+j-1);
    \draw[gray!60,line width=3mm] (ai+3)    -- (7,0);
    \draw[thick] (bi+3) -- (7,3);

  \end{scope}

  \begin{scope}
    \clip (8.7,3) rectangle (11,-.1);
    \draw[thick] (7,3) .. controls (8,.5) and (10,1.3) .. (ai+j);
    \draw[thick] (bi+j-1) -- (8,3);
    \draw[gray!60,line width=3mm] (ai+j-1) -- (8,0);
    \end{scope}

  \draw[thick]
    (ci+j) .. controls (10,.5) and (12,1.3) .. (ai+j+1);

  \begin{scope}
    \clip (12,3) rectangle (13.3,-.1);
    \draw[thick] (ci+j+1) .. controls (12,.5) and (14,1.3) .. (16,0);
    \draw[thick] (bi+j+1)--(14,3);
    \draw[thick] (ai+j+1)--(14,0);
  \end{scope}

  \node[fill=black,circle,inner sep=2pt,label=below:$a_k$]           at (ai)       {};
  \node[fill=black,circle,inner sep=2pt,label=below:$a_{k+1}$]       at (ai+1)     {};
  \node[fill=black,circle,inner sep=2pt,label=below:$a_{k+2}$]       at (ai+2)     {};
  \node[fill=black,circle,inner sep=2pt,label=below:$a_{k+3}$]       at (ai+3)     {};
  \node[fill=black,circle,inner sep=2pt,label=below:$a_{k+l-1}$]     at (ai+j-1)   {};
  \node[fill=black,circle,inner sep=2pt,label=below:$a_{k+l}$]       at (ai+j)     {};
  \node[fill=black,circle,inner sep=2pt,label=below:$a_{k+l+1}$]     at (ai+j+1)   {};

  \node[fill=black,circle,inner sep=2pt,label=-10:$b_k$]            at (bi)       {};
  \node[fill=black,circle,inner sep=2pt,label=-10:$b_{k+1}$]        at (bi+1)     {};
  \node[fill=black,circle,inner sep=2pt,label=-10:$b_{k+2}$]        at (bi+2)     {};
  \node[fill=black,circle,inner sep=2pt,label=-10:$b_{k+3}$]        at (bi+3)     {};
  \node[fill=black,circle,inner sep=2pt,label=-10:$b_{k+l-1}$]      at (bi+j-1)   {};
  \node[fill=black,circle,inner sep=2pt,label=-10:$b_{k+l}$]        at (bi+j)     {};
  \node[fill=black,circle,inner sep=2pt,label=-10:$b_{k+l+1}$]      at (bi+j+1)   {};

  \node[fill=black,circle,inner sep=2pt,label=above:$c_{k+1}$]       at (ci+1)     {};
  \node[fill=black,circle,inner sep=2pt,label=above:$c_{k+2}$]       at (ci+2)     {};
  \node[fill=black,circle,inner sep=2pt,label=above:$c_{k+3}$]      at (ci+3)     {};
  \node[fill=black,circle,inner sep=2pt,label=above:$c_{k+l}$]      at (ci+j)     {};
  \node[fill=black,circle,inner sep=2pt,label=above:$c_{k+l+1}$]    at (ci+j+1)   {};

\end{tikzpicture}
}
    \caption{Type (3,a)}\label{fig:3a}
  \end{subfigure}\hfill
  \begin{subfigure}[t]{0.48\linewidth}
    \centering
    \resizebox{\linewidth}{!}{\begin{tikzpicture}[
  every node/.style={font=\Large},
  every label/.append style={font=\Large},
  baseline=(current bounding box.north) 
]
  \path[use as bounding box] (-1.4,-1.0) rectangle (16.4,4.3);

  \node at (7.5,1)       {\Huge $\cdots$};

  \coordinate (ak)       at ( 0,0);
  \coordinate (ak+1)     at ( 2,0);
  \coordinate (ak+2)     at ( 4,0);
  \coordinate (ak+3)     at ( 6,0);
  \coordinate (ak+j-1)   at ( 9,0);
  \coordinate (ak+j)     at (11,0);
  \coordinate (ak+j+1)   at (13,0);

  \coordinate (bk)       at ( 0,2);
  \coordinate (bk+1)     at ( 2,2);
  \coordinate (bk+2)     at ( 4,2);
  \coordinate (bk+3)     at ( 6,2);
  \coordinate (bk+j-1)   at ( 9,2);
  \coordinate (bk+j)     at (11,2);
  \coordinate (bk+j+1)   at (13,2);

  \coordinate (ck)       at ( -1,3);
  \coordinate (ck+1)     at ( 1,3);
  \coordinate (ck+2)     at ( 3,3);
  \coordinate (ck+3)     at ( 5,3);
  \coordinate (ck+j)     at (10,3);
  \coordinate (ck+j+1)   at (12,3);

  \draw[gray!60,line width=3mm]
    (bk) -- (ck+1)
    (ck+1) .. controls (1,.5) and (3,1.3) .. (ak+2) -- (ak+3)
    (ak+3) -- (bk+3)
    (bk+j-1) -- (ck+j) -- (bk+j);

  \draw[thick]
    (ak)      -- (ak+1) -- (ak+2)
    (ak+j-1)  -- (ak+j)  -- (ak+j+1);

  \draw[thick]
    (ak)      -- (bk)
    (ak+1)    -- (bk+1)
    (ak+2)    -- (bk+2)
    (ak+j-1)  -- (bk+j-1)
    (ak+j)    -- (bk+j)
    (ak+j+1)  -- (bk+j+1);

  \draw[thick]
    (ck+1) -- (bk+1) -- (ck+2) -- (bk+2) -- (ck+3) -- (bk+3)
    (bk+j)  -- (ck+j+1) -- (bk+j+1);

  \begin{scope}
    \clip (-.3,0) rectangle (2.1,3.1);
    \draw[thick]
      (bk) -- (ck)
      (ck) .. controls (-1,.5) and (1,1.3) .. (ak+1)
      (ak) -- (-1,0);
  \end{scope}

  \draw[thick]
    (ck+2) .. controls (3,.5)  and (5,1.3)  .. (ak+3);

  \begin{scope}
    \clip (4,3) rectangle (6.3,-.1);
    \draw[thick]
      (ck+3) .. controls (5,.5) and (7,1.3) .. (ak+j-1);
      \draw[thick](ak+3) -- (7,0);
    \draw[gray!60,line width=3mm] (bk+3) -- (7,3);
  \end{scope}

  \begin{scope}
    \clip (8.7,3) rectangle (11,-.1);
    \draw[thick]
      (7,3) .. controls (8,.5)   and (10,1.3) .. (ak+j)
      (bk+j-1) -- (8,3);
    \draw[thick] (ak+j-1) -- (8,0);
    \draw[gray!60,line width=3mm] (bk+j-1) -- (8,3);

  \end{scope}

  \draw[thick]
    (ck+j) .. controls (10,.5) and (12,1.3) .. (ak+j+1);

  \begin{scope}
    \clip (12,3) rectangle (13.3,-.1);
    \draw[thick]
      (ck+j+1) .. controls (12,.5) and (14,1.3) .. (16,0)
      (bk+j+1) -- (14,3)
      (ak+j+1) -- (14,0);
  \end{scope}

  \node[fill=black,circle,inner sep=2pt,label=below:$a_k$]         at (ak)       {};
  \node[fill=black,circle,inner sep=2pt,label=below:$a_{k+1}$]     at (ak+1)     {};
  \node[fill=black,circle,inner sep=2pt,label=below:$a_{k+2}$]     at (ak+2)     {};
  \node[fill=black,circle,inner sep=2pt,label=below:$a_{k+3}$]     at (ak+3)     {};
  \node[fill=black,circle,inner sep=2pt,label=below:$a_{k+l-1}$]   at (ak+j-1)   {};
  \node[fill=black,circle,inner sep=2pt,label=below:$a_{k+l}$]     at (ak+j)     {};
  \node[fill=black,circle,inner sep=2pt,label=below:$a_{k+l+1}$]   at (ak+j+1)   {};

  \node[fill=black,circle,inner sep=2pt,label=-10:$b_k$]          at (bk)       {};
  \node[fill=black,circle,inner sep=2pt,label=-10:$b_{k+1}$]      at (bk+1)     {};
  \node[fill=black,circle,inner sep=2pt,label=-10:$b_{k+2}$]      at (bk+2)     {};
  \node[fill=black,circle,inner sep=2pt,label=-10:$b_{k+3}$]      at (bk+3)     {};
  \node[fill=black,circle,inner sep=2pt,label=-10:$b_{k+l-1}$]    at (bk+j-1)   {};
  \node[fill=black,circle,inner sep=2pt,label=-10:$b_{k+l}$]      at (bk+j)     {};
  \node[fill=black,circle,inner sep=2pt,label=-10:$b_{k+l+1}$]    at (bk+j+1)   {};

  \node[fill=black,circle,inner sep=2pt,label=above:$c_{k+1}$]    at (ck+1)     {};
  \node[fill=black,circle,inner sep=2pt,label=above:$c_{k+2}$]    at (ck+2)     {};
  \node[fill=black,circle,inner sep=2pt,label=above:$c_{k+3}$]    at (ck+3)     {};
  \node[fill=black,circle,inner sep=2pt,label=above:$c_{k+l}$]    at (ck+j)     {};
  \node[fill=black,circle,inner sep=2pt,label=above:$c_{k+l+1}$]  at (ck+j+1)   {};

\end{tikzpicture}
}
    \caption{Type (3,b)}\label{fig:3b}
  \end{subfigure}

  \caption{Induced paths of length $0, 1,$ and $2 \mod 3$ in $G_m$.}
  \label{fig:allPaths}
\end{figure}
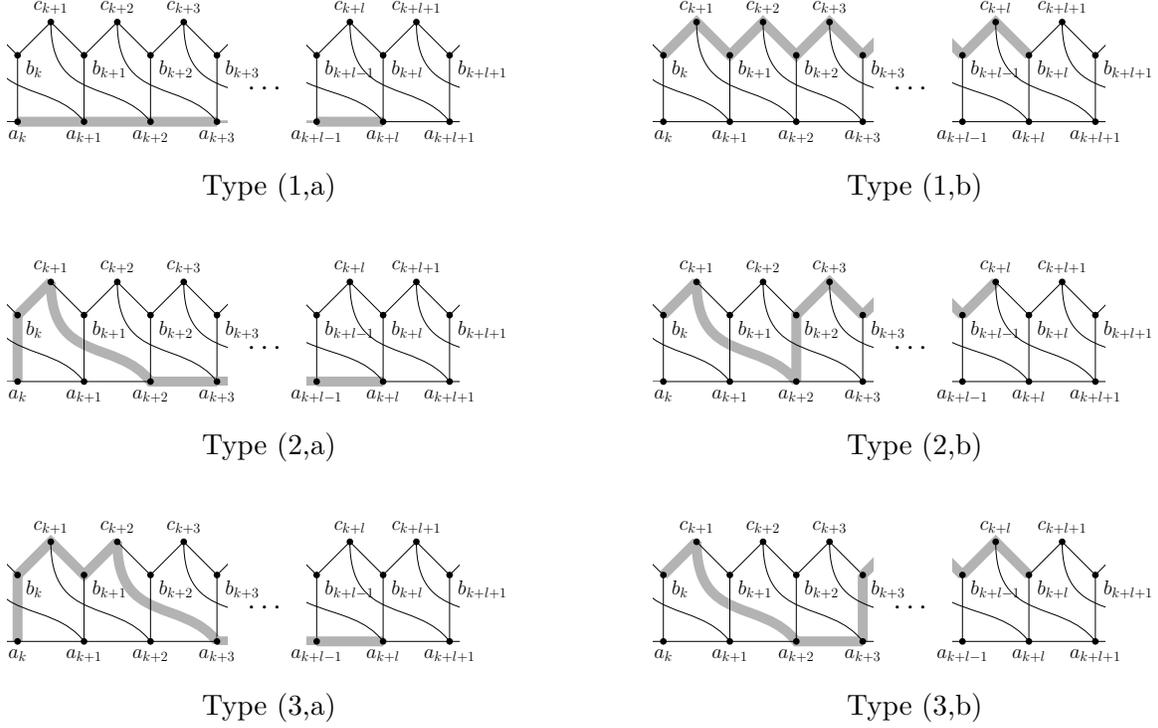

    In the following, all paths are induced paths unless otherwise stated.
    We start by finding paths of different lengths between $a_k$ and $a_{k+l}$, and between $b_k$ and $b_{k+l}$. These paths are a crucial step in finding the desired paths claimed in the statement of the lemma.
    Let $k \in [m]$ and $l \geq 0$ such that $k+l \leq m$.
    A \emph{type--$(1,a)$--path} is a path of length $l$ of the form $(a_k, a_{k+1}, \ldots, a_{k+l})$, where $l \geq 1$. A \emph{type--$(2,a)$--path} is a path of length $l+1$ of the form $(a_k, b_k, c_{k+1}, a_{k+2}, a_{k+3}, \ldots, a_{k+l})$, where $l \geq 2$. A \emph{type--$(3,a)$--path} is a path of length $l+2$ of the form $(a_k, b_k, c_{k+1}, b_{k+1}, c_{k+2}, a_{k+3}, a_{k+4}, \ldots, a_{k+l})$, where $l \geq 3$.  (See Figure~\ref{fig:allPaths}.)
    Observe that these three paths immediately give the desired paths if $x,y \in \{a_i\}_{i \in [m]}$.
    Moreover, if $x=a_i$ and $y = b_j$, we obtain the desired paths as follows. Note that since $\dist(x,y) \geq 3$, we may disregard $j \in \{i-2, i-1, i, i+1\}$.
    If $j \geq i+3$, concatenating a type--$(t,a)$--path for $t \in [3]$ with the edge $a_jb_j$ gives the desired paths.
    Similarly, if $j \leq i-4$, concatenating the path $(b_j,a_j,a_{j+1})$ with a type--$(t,a)$--path for $t \in [3]$ gives the desired paths.
    For the remaining two cases, we refer to Table \ref{tab:paths_ai}.

    If $x=a_i$ and $y=c_j$, we obtain the desired paths as follows. Again, since $\dist(x,y) \geq 3$, we may disregard $j \in \{i-2, i-1, i, i+1\}$.
    If $j \geq i+3$, concatenating a type--$(t,a)$--path for $t \in [3]$ with the path $(a_j,b_j,c_j)$ gives the desired paths.
    Similarly, if $j \leq i-4$, concatenating the edge $c_ja_{j+1}$ with a type--$(t,a)$--path for $t \in [3]$ gives the desired paths.
    For the remaining two cases, we refer to Table \ref{tab:paths_ai}.
    \begin{table}[H]
        \centering
        \begin{tabular}{c|c||c|c|c}
            $x$ & $y$   & $0 \mod 3$    & $1 \mod 3$    & $2 \mod 3$ \\
            \hline\hline
            $a_i$   & $b_{i+2}$ & $(a_i,a_{i+1},a_{i+2},b_{i+2})$    & $(a_i,b_i,c_{i+1},a_{i+2},b_{i+2})$     & $(a_i,b_i,c_{i+1},b_{i+1},c_{i+2},b_{i+2})$ \\
            \hline
            $a_i$   & $b_{i-3}$ & $(b_{i-3},c_{i-2},a_{i-1},a_i)$    & $(b_{i-3},a_{i-3},a_{i-2},a_{i-1},a_i)$ & $(b_{i-3},a_{i-3},a_{i-2},b_{i-2},c_{i-1},a_i)$\\
            \hline
            $a_i$   & $c_{i+2}$ & $(a_i,a_{i+1},b_{i+1},c_{i+2})$    & $(a_i,b_i,c_{i+1},b_{i+1},c_{i+2})$     & $(a_i,b_i,c_{i+1},a_{i+2},b_{i+2},c_{i+2})$\\
            \hline
            $a_i$   & $c_{i-3}$ & $(c_{i-3},a_{i-2},a_{i-1},a_i)$    & $(c_{i-3},b_{i-3},c_{i-2},a_{i-1},a_i)$ & $(c_{i-3},b_{i-3},a_{i-3},a_{i-2},a_{i-1},a_i)$
        \end{tabular}
        \caption{Short paths: $a_i$}
        \label{tab:paths_ai}
    \end{table}

    Similarly, we can define a \emph{type--$(1,b)$--path} as a path of length $2l$ of the form $(b_k, c_{k+1},$ $b_{k+1}, \ldots, c_{k+l}, b_{k+l})$, where $l \geq 1$.
    A \emph{type--$(2,b)$--path} is a path of length $2l-1$ of the form $(b_k, c_{k+1}, a_{k+2}, b_{k+2}, c_{k+3}, b_{k+3}, \ldots, c_{k+l},$ $b_{k+l})$, where $l \geq 2$.
    A \emph{type--$(3,b)$--path} is a path of length $2l-2$ of the form $(b_k, c_{k+1}, a_{k+2}, a_{k+3}, b_{k+3}, c_{k+4}, b_{k+4}, \ldots,$ $c_{k+l}, b_{k+l})$, where $l \geq 3$.  (See Figure~\ref{fig:allPaths}.)
    These immediately give the desired paths if $x=b_i, y=b_j$, and $|i-j| \geq 3$. As $\dist(x,y) \geq 3$, we may disregard pairs $i,j \in [m]$ with $|i-j| \leq 1$. If $|i-j|=2$, assume without loss of generality that $i<j$, and we refer to Table \ref{tab:paths_bi}. 
    
    If $x=b_i$ and $y=c_j$, we obtain the desired paths as follows. Note that since $\dist(x,y) \geq 3$, we may disregard $j \in \{i-1, i, i+1\}$.
    If $j \geq i+4$, concatenating a type--$(t,b)$--path for $t \in [3]$ with the edge $b_jc_j$ gives the desired paths. If $j \leq i-3$, concatenating the edge $c_jb_j$ with a type--$(t,b)$--path for $t \in [3]$ gives the desired paths. For the remaining three cases, we refer to Table \ref{tab:paths_bi}.
    \begin{table}[H]
        \centering
        \scriptsize
        \begin{tabular}{c|c||c|c|c}
            $x$     & $y$                                                       & $0 \mod 3$                                    & $1 \mod 3$                            & $2 \mod 3$ \\
            \hline \hline
            $b_i$   & $b_{i+2}$                                                 & $(b_i,c_{i+1},a_{i+2},b_{i+2})$                    & $(b_i,c_{i+1},b_{i+1},c_{i+2},b_{i+2})$     & $(b_i,a_i,a_{i+1},b_{i+1},c_{i+2},b_{i+2})$\\
            \hline
            $b_i$   & $c_{i+3}$                                                 & $(b_i,a_i,a_{i+1},a_{i+2},a_{i+3},b_{i+3},c_{i+3})$   & $(b_i,c_{i+1},a_{i+2},b_{i+2},c_{i+3})$     & $(b_i,c_{i+1},b_{i+1},c_{i+2},b_{i+2},c_{i+3})$\\
            \hline
            $b_i$   & $c_{i+2}$                                                 & $(b_i,c_{i+1},b_{i+2},c_{i+2})$                    & $(b_i,a_i,a_{i+1},b_{i+1},c_{i+2})$         & $(b_i,a_i,a_{i+1},a_{i+2},b_{i+2},c_{i+2})$\\
            \hline
            $b_i$   & $c_{i-2}$                                                 & $(c_{i-2},a_{i-1},a_{i},b_i)$                    & $(c_{i-2},b_{i-2},c_{i-1},a_i,b_i)$         & $(c_{i-2},b_{i-2},c_{i-1},b_{i-1},c_i,b_i)$
        \end{tabular}
        \caption{Short paths: $b_i$}
        \label{tab:paths_bi}
    \end{table}  

    Finally, if $x=c_i, y=c_j$, assume without loss of generality that $i < j$. As $\dist(x,y) \geq 3$, we may disregard pairs $i,j \in [m]$ with $j \leq i+1$. If $j \geq i+4$, concatenating the edge $c_ib_i$, a type--$(t,b)$--path for $t \in [3]$, and the edge $b_{j-1}c_j$ gives the desired paths. For the remaining two cases, we refer to Table \ref{tab:paths_ci}.
    \begin{table}[H]
        \centering
        \scriptsize
        \begin{tabular}{c|c||c|c|c}
            $x$     & $y$       & $0 \mod 3$                                    & $1 \mod 3$                        & $2 \mod 3$ \\
            \hline\hline
            $c_i$   & $c_{i+3}$ & $(c_i,b_i,c_{i+1},b_{i+1},c_{i+2},b_{i+2},c_{i+3})$   & $(c_i,a_{i+1},a_{i+2},b_{i+2},c_{i+3})$ & $(c_i,a_{i+1},a_{i+2},a_{i+3},b_{i+3},c_{i+3})$ \\
            \hline
            $c_i$   & $c_{i+2}$ & $(c_i,a_{i+1},b_{i+1},c_{i+2})$                    & $(c_i,b_i,c_{i+1},b_{i+1},c_{i+2})$     & $(c_i,b_i,c_{i+1},a_{i+2},b_{i+2},c_{i+2})$
        \end{tabular}
        \caption{Short paths: $c_i$}
        \label{tab:paths_ci}
    \end{table}
\end{proof}

The recursive construction of $G_m$ allows us to recursively describe the independent sets of $G_m$. We can moreover specify when such independent sets are maximal.

\begin{proposition} \label{prop:indep_gm}
    Every independent set $A$ of $G_m$ has one of the following forms:
    \begin{enumerate}
        \item $A \cap \{a_m, b_m, c_m\}=\emptyset$ and $A$ is an independent set in $G_k$ for some $k < m$.
        \item $A \cap \{a_m, b_m, c_m\}=\{b_m\}$ and $A \setminus \{b_m\}$ is an independent set in $G_k$ for some $k < m$.
        \item $A \cap \{a_m, b_m, c_m\}=\{a_m, c_m\}$ and $A \setminus \{a_m, c_m\}$ is an independent set in $G_k$ for some $k \leq m-2$. In particular, $A \cap \{a_{m-1}, b_{m-1}, c_{m-1}\}=\emptyset$.
        \item $A \cap \{a_m, b_m, c_m\}=\{a_m\}$ and $A \setminus \{a_m\}$ is an independent set in $G_k$ for some $k < m$ disjoint from $\{a_{m-1}, c_{m-1}\}$.
        \item $A \cap \{a_m, b_m, c_m\}=\{c_m\}$ and $A \setminus \{c_m\}$ is an independent set in $G_k$ for some $k < m$ disjoint from $\{b_{m-1}\}$.
    \end{enumerate}
    Moreover, the independent sets from (2), (4), and (5) are maximal if $A \setminus \{b_m\}$, $A \setminus \{a_m\}$, and $A \setminus \{c_m\}$ are maximal independent sets in $G_{m-1}$, respectively. Finally, the independent set from (3) is maximal if $A \setminus \{a_m, c_m\}$ is a maximal independent set in $G_{m-2}$.
\end{proposition}
\begin{proof}
    The independent sets (1) - (5) immediately follow from the construction of $G_m$. To obtain the maximality statements, recall that a maximal independent set in $G_m$ has cardinality $m$.
\end{proof}

\begin{definition}\label{def:edgesubdivision} 
Let $\Delta$ be a simplicial complex and $\{x,y\}$ an edge of $\Delta$.  The \emph{edge subdivision} of $\Delta$ at $\{x,y\}$ is the simplicial complex $\Delta'$ with vertex set $V(\Delta)\sqcup\{a\}$ and faces:
\begin{itemize}
    \item $\sigma$, for $\sigma\in \Delta$ such that $\{x,y\}\not\subset \sigma$
\item $\sigma \setminus \{x\}\cup\{a\}$, for $\sigma \in \Delta$ such that $\{x,y\} \subset \sigma$
\item $\sigma \setminus \{y\}\cup\{a\}$, for $\sigma \in \Delta$ such that $\{x,y\} \subset \sigma$
\end{itemize}
\end{definition}

\begin{remark}\label{rem:edgesubdivision}
Edge subdivisions preserve flagness of the complex (\cite[Lemma 5.2]{lutznevo16}). We will focus on performing edge subdivisions on $\Delta$ when $\Delta = \Ind(G)$ for some graph $G$. We see how an edge subdivision of $\Delta$ comes from an operation on $G$. 
Observe that $G$ is the complement of the $1$-skeleton of $\Delta$. 
Therefore, using the notation from the definition above, $\Delta' = \Ind(G')$ where $G'$ is a graph with vertex set $V(G) \sqcup \{a\}$ and edges
\begin{itemize}
    \item $xy$
    \item $e$, for every $e \in E(G)$
    \item $az$, for every $z \in N_G(x) \cup N_G(y)$
\end{itemize}   
In \cref{fig:Ind(G_3)}, the right-hand-side sphere is obtained by performing an edge subdivision at the edge $a_1a_2$ on the left-hand-side sphere with a new vertex $c_2$. On the left-hand-side graph, $N(a_1) \cup N(a_2) = \{b_1, b_2, a_3\}$. Therefore, the edge subdivision is reflected on the graph as adding a vertex $c_3$ and edges $a_1a_2, c_2b_1, c_2b_2, c_2a_3$. 
\end{remark}

A special polytope that will play an important role in the following results in this section is the $m$-dimensional \emph{crosspolytope}, the convex hull of $\{\vec e_i,-\vec e_i : i=1,\ldots,m\}$ in $\R^m$.
Moreover, any polytope combinatorially equivalent to a crosspolytope is called a crosspolytope.  The independence complex of the graph of $m$ disjoint edges ($mK_2$) is the boundary of a crosspolytope.

\begin{theorem}\label{thm:crosspolysubdivide}
 For each $m\ge 1$   the simplicial complex $\Ind(G_m)$ is obtained by a 
sequence of $m-1$ edge subdivisions of the boundary of the
$m$-dimensional crosspolytope.
\end{theorem}

\begin{proof}
Fix $m$ and let $Q^m_1$ be the boundary of the $m$-dimensional crosspolytope.
$Q^m_1$ is the independence complex of 
$mK_2$, which is $mG_1$;
call this graph $S_1$.
Let the vertex set of $S_1$ be $\{b_i: 1\le i\le m\}\cup \{c_i: 1\le i\le m\}$, with edge set $\{b_ic_i: 1\le i\le m\}$.
Now we perform a sequence of edge subdivisions on $Q^m_1$,  producing a sequence of flag complexes, $Q^m_i$, and of graphs $S_i$, the complements of the 1-skeletons.  (These are the graphs for which $Q^m_i$ are the independence complexes.) 
Subdivide edge $b_1c_2$ of $Q^m_1$ with new vertex $a_2$.
Call the new complex $Q^m_2$, and the complement of its 1-skeleton $S_2$.
The graph $S_2$ has edges $b_ic_i$ and $b_1c_2$, $a_2c_1$, and $a_2b_2$.
Thus $S_2$ is a pentagon along with $m-2$ disjoint edges, that is, $S_2$ is
the disjoint union of $G_2$ and $m-2$ copies of $G_1$.
Suppose we have formed $Q^m_i$ and $S_i$, $2\le i\le m-2$, and we know that $S_i$ is the disjoint union of $G_i$ and $m-i$ copies of $G_1$.  Then subdivide the edge
$b_ic_{i+1}$ with new vertex $a_{i+1}$.
Call the new complex $Q^m_{i+1}$ and the complement of its 1-skeleton $S_{i+1}$.
$S_{i+1}$ has the edges of $S_i$ (including $b_{i+1}c_{i+1}$), and edges
$b_ic_{i+1}$ and  $a_{i+1}c_i$, $a_{i+1}b_{i+1}$, and $a_ia_{i+1}$.
Thus $S_{i+1}$ is the disjoint union of $G_{i+1}$ and $m-i-1$ copies of $G_1$.
Finally, when $i=m-1$, $S_m=G_m$.
So $G_m$ is the complement of the 1-skeleton of $Q^m_m$.
Thus, $Q^m_m$, the simplicial complex resulting from subdividing $m-1$ edges of the boundary of the $m$-dimensional crosspolytope, is the independence complex of $G_m$.
\end{proof} 

In order to prove Theorem \ref{t:onlygmintro}, we establish the two corollaries below, which follow from \cite[Proposition 2.4.2]{gal2005} and \cite[Theorem 2.7]{provan-billera}, respectively.

\begin{corollary}
For $m\ge 1$, $\Ind(G_m)$ is the boundary of an $m$-dimensional simplicial 
polytope. 
\end{corollary}

\begin{corollary}
For $m\ge 1$, $\Ind(G_m)$ is vertex decomposable.
\end{corollary}

Theorem \ref{t:onlygmintro} now quickly follows from Theorem \ref{t:trungcharacterization} and the two previous corollaries.

\thmB*

\begin{proof}
Let $G$ be a planar ternary graph such that $\Ind(G)$ is a homology sphere.  Then $G$ is Gorenstein, so by \cref{t:trungcharacterization}, $G$ is the disjoint union of $G_m$ graphs.  Thus $\Ind(G)$ is the join of the independence complexes $\Ind(G_m)$, each of which is the boundary of a simplicial polytope and is vertex decomposable.  These properties are preserved under the join operation, so $\Ind(G)$ is the boundary of a simplicial polytope and is vertex decomposable. 
\end{proof}


\section{An induced subgraph of the Lutz--Nevo graph} \label{sec:FlipGraph}

\cref{thm:crosspolysubdivide} demonstrates how to obtain $\Ind(G_m)$ from the boundary of the crosspolytope by edge subdivisions. As a result, independence complexes of disjoint unions of $G_m$'s are polytopal simplicial spheres. Given two such spheres $\Ind(G)$ and $\Ind(G')$ of the same dimension, we now show how to obtain one from the other by edge subdivisions and contractions. For example, \cref{fig:Ind(G_3)} illustrates an edge subdivision that turns $\Ind(G_1 \sqcup G_2)$ into $\Ind(G_3)$. These observations allow us 
to show that the subgraph of the Lutz--Nevo graph induced on the set of vertices representing these spheres is isomorphic to the partition refinement graph.

\begin{theorem}[{\cite[Theorem 1.2]{lutznevo16}}]\label{t:connected_flip}
    Two flag simplicial complexes are PL homeomorphic if and only if they can be connected by a sequence of edge contractions and subdivisions such that every complex in the sequence is flag.
\end{theorem}

In view of~\cref{t:connected_flip} the goal of this section is to study the following graph.

\begin{definition}
    Let $V_d$ be the set of all isomorphism classes of $d$-dimensional flag PL spheres for some fixed $d$, and $E_d$ be the collection of pairs $\{\Delta_1, \Delta_2\}$ of elements of $V_d$ such that $\{\Delta_1, \Delta_2\} \in E_d$ if and only if $\Delta_2$ can be obtained from $\Delta_1$ by either contracting an edge, or subdividing an edge. For an integer $d \geq 0$, the \emph{Lutz--Nevo graph} $T_d$ is the graph with vertex set $V_d$ and edge set $E_d$.    
\end{definition}

\begin{lemma}[\cite{lutznevo16}]
    Let $\Delta$ be a flag complex. The complex $\Delta'$ obtained by contracting an edge $e$ of $\Delta$ is flag if and only if $e$ is not part of an induced square of $\Delta$.
\end{lemma}

We are interested in the subgraph associated with the planar ternary graphs.
\begin{definition}
    Denote by $H_{n-1}$ the induced subgraph of $T_{n-1}$ where the vertices are given by the (isomorphism types) of  independence complexes of planar ternary graphs with independence number $n$. 
\end{definition}

We will show this subgraph is isomorphic to a graph arising from partitions.

A \emph{partition} $\lambda = (m_1, \dots, m_s)$ of $n$ is an ordered nonincreasing list of integers such that $m_1 + \dots + m_s = n$.

\begin{definition}[{\bf 
Partition refinement graph}]
     A partition $\lambda$ of $n$ \textit{refines} a partition $\mu$ of $n$ if each part of $\mu$ is the sum of parts of $\lambda$. The \textit{partition refinement graph}, $\Pp_n$, is the graph of the Hasse diagram of the refinement poset on the integer partitions of $n$.
\end{definition}

We will prove that $\Pp_n$ is isomorphic to $H_{n-1}$.

\begin{lemma}
    There is a bijection between $V(H_{n-1})$ and $V(\Pp_n)$.    
\end{lemma}

\begin{proof}
    If the independence complex of a planar ternary graph $G$ is an $(n-1)$-dimensional sphere, then by~\cite{trung2018characterization} 
    $$
        G = G_{m_1} \sqcup G_{m_2} \sqcup \dots \sqcup G_{m_s}
    $$
    where $m_1 + \dots + m_s = n$. Reordering the $m_i$ does not change the isomorphism class. Hence the map from $V(H_{n-1})$ to $V(\Pp_n)$ given by 
    $$
        \Ind(G_{m_1} \sqcup G_{m_2} \sqcup \dots \sqcup G_{m_s}) \mapsto (m_1, \dots, m_s)
    $$ is a bijection.
\end{proof}

\begin{lemma}\label{lem:edge_in_subdivision}
    Let $\Ind(G)$ be a vertex of $H_{n - 1}$, and $xy$ an edge of $\Ind(G)$ where $x$ and $y$ are vertices of the same connected component in $G$. Then the sphere given by subdividing the edge $xy$ is not a vertex of $H_{n-1}$.
\end{lemma}

\begin{proof}
    Subdividing the edge $xy$ of $\Ind(G_m)$ adds the edge $xy$ to $G_m$ (among other edges\, /\,vertices). There are two cases:

    \begin{enumerate}
        \item If the distance between $x$ and $y$ is two, then the new graph has a triangle, and hence is not ternary.
        \item If the distance between $x$ and $y$ is at least three, by~\cref{l:paths_mod3} there exists an induced path between $x$ and $y$ of length $2$ mod $3$. In particular, this path together with the new edge $xy$ gives us an induced cycle of length $3k$ and hence the graph is not ternary.
    \end{enumerate}
\end{proof}

\begin{lemma}\label{l:noextraedges}
Let $G = G_{m_1} \sqcup G_{m_2} \sqcup \dots \sqcup G_{m_s}$ and $
\Delta = \Ind(G)$. 
Suppose $xy$ is an edge of $\Delta$, where $x$ is a vertex of $G_{m_i}$ and
$y$ is a vertex of $G_{m_j}$ ($i\ne j$).  Let $\Delta'$ be obtained from $\Delta$ by an
edge subdivision of $xy$, and let $G'$ be the graph with $\Ind(G')=\Delta'$.
Then $G'$ has a connected component with subgraphs $G_{m_i}$ and $G_{m_j}$,
and this component has independence number $m_i+m_j$.
\end{lemma}

\begin{proof}
Let $\Delta'$ be as given in the statement of the lemma, with new vertex
$a$ created by the edge subdivision.  
$\Delta'$ is a flag complex, so it is
the independence complex of some graph $G'$.  The faces of $\Delta'$ are those 
faces of $\Delta$ not containing both $x$ and $y$, and the faces 
$\sigma \setminus \{x\} \cup \{a\}$ and
$\sigma \setminus \{y\} \cup \{a\}$ for each $\sigma\in\Delta$ containing 
$\{x,y\}$.  
The edges of $G'$ are the nonedges of $\Delta'$; these are the nonedges of 
$\Delta$, the pair $xy$ and some edges containing $a$.  

Note that if $u$ is a vertex
of $G$ not in $V(G_{m_i}\cup G_{m_j})$, then $ux\in \Delta$, so
$ua\in\Delta'$, so $ua$ is not an edge of $G'$. Thus, the
subgraph $G''$ induced by the vertex set $V(G_{m_i}\cup G_{m_j}) \cup \{a\}$
is a connected component of $G'$ containing as subgraphs $G_{m_i}$ and 
$G_{m_j}$.  If $A_i$ is a maximal independent set of $G_{m_i}$ containing $x$,
and $A_j$ is a maximal independent set of $G_{m_j}$ containing $y$, then
$A_i\cup A_j \setminus $\{x\}$ \cup \{a\}$ is an independent set of size 
$m_i+m_j$ in $G'$.  Moreover, this set must be a maximal independent set in
$G''$, because every other vertex of $G_{m_i}\cup G_{m_j}$ is adjacent to
some vertex of $A_i\cup A_j$.
\end{proof}

\thmC*

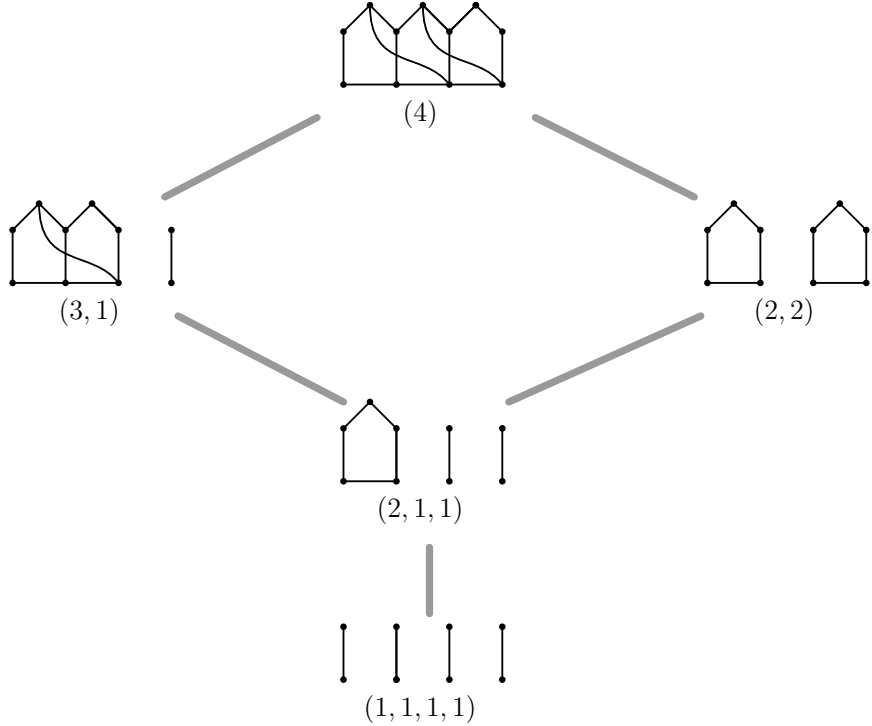
\begin{figure}[htbp]
    \centering
    \resizebox{.7\textwidth}{!}{\begin{tikzpicture}[
  vertex/.style={fill=black, circle, inner sep=1pt},
  every label/.append style={font=\small, label distance=-2pt},
  thick
]

\begin{scope}[shift={(-.5,0)}, scale=0.4] 
\coordinate (ai-2) at (-4,0);
\coordinate (ai-1) at (-2,0);
\coordinate (ai)   at (0,0);
\coordinate (ai+1) at (2,0);
\coordinate (ai+2) at (4,0);

\coordinate (bi-2) at (-4,2);
\coordinate (bi-1) at (-2,2);
\coordinate (bi)   at (0,2);
\coordinate (bi+1) at (2,2);
\coordinate (bi+2) at (4,2);

\coordinate (ci)   at (-1,3);
\coordinate (ci+1) at (1,3);
\coordinate (ci+2) at (3,3);

\draw (ai-1)--(bi-1) (ai)--(bi) (ai+1)--(bi+1) (ai+2)--(bi+2);
\draw (bi-1)--(ci)--(bi)--(ci+1)--(bi+1)--(ci+2)--(bi+2);
\draw (ci)   .. controls (-1,.5) and (1,1.3) .. (ai+1);
\draw (ci+1) .. controls ( 1,.5) and (3,1.3) .. (ai+2);
\draw (ci+1)--(bi+1)--(ci+2);
\draw (ai-1)--(ai+2);

\node[vertex] at (ai-1) {};
\node[vertex] at (ai)   {};
\node[vertex] at (ai+1) {};
\node[vertex] at (ai+2) {};

\node[vertex] at (bi-1) {};
\node[vertex]   at (bi)   {};
\node[vertex]   at (bi+1) {};
\node[vertex]   at (bi+2) {};

\node[vertex]  at (ci)   {};
\node[vertex]  at (ci+1) {};
\node[vertex]  at (ci+2) {};

\node at (0.9,-1.1) {$(4)$};
\end{scope}

\begin{scope}[shift={(-5.5,-3)}, scale=0.4]
\coordinate (ai-2) at (-4,0);
\coordinate (ai-1) at (-2,0);
\coordinate (ai)   at (0,0);
\coordinate (ai+1) at (2,0);
\coordinate (ai+2) at (4,0);

\coordinate (bi-2) at (-4,2);
\coordinate (bi-1) at (-2,2);
\coordinate (bi)   at (0,2);
\coordinate (bi+1) at (2,2);
\coordinate (bi+2) at (4,2);

\coordinate (ci) at (-1,3);
\coordinate (ci+1) at (1,3);

\draw (ai-1)--(bi-1) (ai)--(bi) (ai+1)--(bi+1) (ai+2)--(bi+2);
\draw (bi-1)--(ci)--(bi)--(ci+1)--(bi+1);
\draw (ci) .. controls (-1,.5) and (1,1.3) .. (ai+1);
\draw (ci+1)--(bi+1);
\draw (ai-1)--(ai+1);

\node[vertex] at (ai-1) {};
\node[vertex] at (ai)   {};
\node[vertex] at (ai+1) {};
\node[vertex] at (ai+2) {};

\node[vertex] at (bi-1) {};
\node[vertex]   at (bi)   {};
\node[vertex] at (bi+1) {};
\node[vertex]   at (bi+2) {};

\node[vertex] at (ci) {};
\node[vertex]  at (ci+1) {};
\node[vertex]  at (ci+2) {};

\node at (0.9,-1.1) {$(3,1)$};

\end{scope}

\begin{scope}[shift={(5,-3)}, scale=0.4]
\coordinate (ai-2) at (-4,0);
\coordinate (ai-1) at (-2,0);
\coordinate (ai)   at (0,0);
\coordinate (ai+1) at (2,0);
\coordinate (ai+2) at (4,0);

\coordinate (bi-2) at (-4,2);
\coordinate (bi-1) at (-2,2);
\coordinate (bi)   at (0,2);
\coordinate (bi+1) at (2,2);
\coordinate (bi+2) at (4,2);

\coordinate (ci)   at (-1,3);
\coordinate (ci+2) at (3,3);

\draw (ai-1)--(bi-1) (ai)--(bi) (ai+1)--(bi+1) (ai+2)--(bi+2);
\draw (bi-1)--(ci)--(bi) (ci+2)--(bi+2);
\draw (bi+1)--(ci+2);
\draw (ai+1)--(ai+2);
\draw (ai-1)--(ai);

\node[vertex] at (ai-1) {};
\node[vertex] at (ai)   {};
\node[vertex] at (ai+1) {};
\node[vertex] at (ai+2) {};

\node[vertex] at (bi-1) {};
\node[vertex]   at (bi)   {};
\node[vertex] at (bi+1) {};
\node[vertex]   at (bi+2) {};

\node[vertex]  at (ci)   {};
\node[vertex]  at (ci+2) {};

\node at (0.9,-1.1) {$(2,2)$};

\end{scope}

\begin{scope}[shift={(-.5,-6)}, scale=0.4]
\coordinate (ai-2) at (-4,0);
\coordinate (ai-1) at (-2,0);
\coordinate (ai)   at (0,0);
\coordinate (ai+1) at (2,0);
\coordinate (ai+2) at (4,0);

\coordinate (bi-2) at (-4,2);
\coordinate (bi-1) at (-2,2);
\coordinate (bi)   at (0,2);
\coordinate (bi+1) at (2,2);
\coordinate (bi+2) at (4,2);

\coordinate (ci) at (-1,3);

\draw (ai-1)--(bi-1) (ai)--(bi) (ai+1)--(bi+1) (ai+2)--(bi+2);
\draw (ai)--(bi);
\draw (bi-1)--(ci)--(bi);
\draw (ai-1)--(ai);

\node[vertex] at (ai-1) {};
\node[vertex] at (ai)   {};
\node[vertex] at (ai+1) {};
\node[vertex] at (ai+2) {};

\node[vertex]   at (bi-1) {};
\node[vertex]   at (bi)   {};
\node[vertex] at (bi+1) {};
\node[vertex]   at (bi+2) {};

\node[vertex]  at (ci) {};

\node at (0.9,-1.1) {$(2,1,1)$};

\end{scope}

\begin{scope}[shift={(-.5,-9)}, scale=0.4]
\coordinate (ai-2) at (-4,0);
\coordinate (ai-1) at (-2,0);
\coordinate (ai)   at (0,0);
\coordinate (ai+1) at (2,0);
\coordinate (ai+2) at (4,0);

\coordinate (bi-2) at (-4,2);
\coordinate (bi-1) at (-2,2);
\coordinate (bi)   at (0,2);
\coordinate (bi+1) at (2,2);
\coordinate (bi+2) at (4,2);

\draw (ai-1)--(bi-1) (ai)--(bi) (ai+1)--(bi+1) (ai+2)--(bi+2);
\draw (ai)--(bi);

\node[vertex] at (ai-1) {};
\node[vertex] at (ai)   {};
\node[vertex] at (ai)   {};
\node[vertex] at (ai+1) {};
\node[vertex] at (ai+2) {};

\node[vertex]   at (bi-1) {};
\node[vertex]   at (bi)   {};
\node[vertex]   at (bi+1) {};
\node[vertex]   at (bi+2) {};

\node at (0.9,-1.1) {$(1,1,1,1)$};

\end{scope}

\draw[line width =1.1mm, gray!80, line cap=round]
    (-1.7,-.5) -- (-4,-1.7)
    (1.6,-.5)  -- (4,-1.7)
    (-3.8,-3.5)  -- (-1.3,-4.8)
    (4.1,-3.5)   -- (1.2,-4.8)
    (0,-7)  -- (0,-8);

\end{tikzpicture}
}
    \caption{A visualization of~\Cref{thm:partitionrefine} for $n = 4$}
    \label{fig:Poset}
    \end{figure}

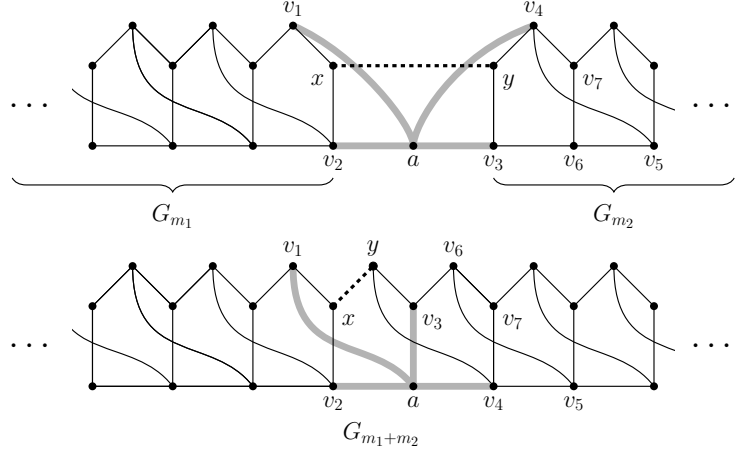
\begin{figure}[H]
    \centering
    \resizebox{.6\textwidth}{!}{

\begin{tikzpicture}[font=\Large]
  \begin{scope}
    \node at (-1.5,1) {\Huge $\cdots$};
    \node at (15.5,1) {\Huge $\cdots$};

    \coordinate (ai) at (0,0);
    \coordinate (ai+1) at (2,0);
    \coordinate (ai+2) at (4,0);
    \coordinate (ai+3) at (6,0);
    \coordinate (ai+j-1) at (8,0);
    \coordinate (ai+j) at (10,0);
    \coordinate (ai+j+1) at (12,0);
    \coordinate (ai+j+2) at (14,0);

    \coordinate (bi) at (0,2);
    \coordinate (bi+1) at (2,2);
    \coordinate (bi+2) at (4,2);
    \coordinate (bi+3) at (6,2);
    \coordinate (bi+j-1) at (8,2);
    \coordinate (bi+j) at (10,2);
    \coordinate (bi+j+1) at (12,2);
    \coordinate (bi+j+2) at (14,2);
    
    \node at (-1,3) (ci) {};
    \coordinate (ci+1) at (1,3);
    \coordinate (ci+2) at (3,3);
    \coordinate (ci+3) at (5,3);
    \coordinate (ci+4) at (7,3);
    \coordinate (ci+j) at (9,3);
    \coordinate (ci+j+1) at (11,3);
    \coordinate (ci+j+2) at (13,3);

        \draw[thick][gray!60, line width=1.7mm] (ci+3) .. controls (6.6,2.3) and (8,.5) .. (ai+j-1);
    \draw[thick] (ai+1) -- (ai+2) -- (ai+3);
    \draw[thick] (ai+j-1) -- (ai+j) -- (ai+j+1);
    \draw[thick] (ai) -- (bi);
    \draw[thick] (ai+1) -- (bi+1);
    \draw[thick] (ai+2) -- (bi+2);
    \draw[thick] (ai+3) -- (bi+3);
    \draw[thick][gray!60, line width=1.7mm] (ai+j-1) .. controls (8,.5) and (9,2.3) .. (ci+j+1);
    \draw[thick] (ai+j) -- (bi+j);
    \draw[thick] (ai+j+1) -- (bi+j+1);
    \draw[thick] (ai) -- (bi) -- (ci+1) -- (bi+1) -- (ci+2) -- (bi+2) -- (ci+3) -- (bi+3);
    \draw[thick][dashed, line width=.8mm] (bi+3) -- (bi+j);
    \draw[thick] (bi+j) -- (ci+j+1) -- (bi+j+1);
    \draw[thick] (ci+1) .. controls (1,.5) and (3,1.3) .. (ai+2);
    \draw[thick] (ci+1) .. controls (1,.5) and (3,1.3) .. (ai+2);
    \draw[thick] (ci+2) .. controls (3,.5) and (5,1.3) .. (ai+3);

    \draw[thick] (ci+1) -- (bi+1) -- (ci+2);
    \draw[thick] (ai) -- (ai+1);
    \draw[thick][gray!60, line width=1.7mm] (ai+3) -- (ai+j);
    \draw[thick] (bi+j) -- (ai+j);
    \draw[thick] (bi+j+1) -- (ci+j+2) -- (bi+j+2) -- (ai+j+2) -- (ai+j+1);

      \draw[thick] (ci+j+1) .. controls (11,.5) and (13,1.3) .. (14,0);

    \node[fill=black, circle, inner sep=2pt] at (ai) {};
    \node[fill=black, circle, inner sep=2pt] at (ai+1) {};
    \node[fill=black, circle, inner sep=2pt] at (ai+2) {};
    \node[fill=black, circle, inner sep=2pt, label=-90:$v_2$] at (ai+3) {};
    \node[fill=black, circle, inner sep=2pt, label=-90:$a$] at (ai+j-1) {};
    \node[fill=black, circle, inner sep=2pt,
  label=-90:$v_3$] at (ai+j) {};
    \node[fill=black, circle, inner sep=2pt,label=-90:$v_6$] at (ai+j+1) {};
    \node[fill=black, circle, inner sep=2pt, label=-90:$v_5$] at (ai+j+2) {};

    \node[fill=black, circle, inner sep=2pt] at (bi) {};
    \node[fill=black, circle, inner sep=2pt] at (bi+1) {};
    \node[fill=black, circle, inner sep=2pt] at (bi+2) {};
    \node[fill=black, circle, inner sep=2pt,
  label=225:$x$] at (bi+3) {};
    
    \node[fill=black, circle, inner sep=2pt,
  label=-45:$y$] at (bi+j) {};
    \node[fill=black, circle, inner sep=2pt, label=-45:$v_7$] at (bi+j+1) {};
    \node[fill=black, circle, inner sep=2pt] at (bi+j+2) {};
    
    \node[fill=black, circle, inner sep=2pt] at (ci+1) {};
    \node[fill=black, circle, inner sep=2pt] at (ci+2) {};
    \node[fill=black, circle, inner sep=2pt, label=90:$v_1$] at (ci+3) {};
    
    \node[fill=black, circle, inner sep=2pt,
  label=90:$v_4$] at (ci+j+1) {};
    \node[fill=black, circle, inner sep=2pt] at (ci+j+2) {};
  \end{scope}

      \begin{scope}
    \clip (-.5,3) rectangle (3.3,-.1);
    \draw[thick] (ci) .. controls (-1,.5) and (1,1.3) .. (2,0);
  \end{scope} 
  
  \begin{scope}
    \clip (12,3) rectangle (14.5,-.1);
    \draw[thick] (ci+j+2) .. controls (13,.5) and (15,1.3) .. (16,0);
  \end{scope}

\draw[thick][decorate,decoration={brace,mirror,amplitude=10pt}]
    (-2,-.8) -- (6,-.8) node[midway,yshift=-27pt] {\Large $G_{m_1}$};

\draw[thick][decorate,decoration={brace,mirror,amplitude=10pt}]
    (10,-.8) -- (16,-.8) node[midway,yshift=-27pt] {\Large $G_{m_2}$};

  \begin{scope}[shift={(0,-6)}]
    \node at (-1.5,1) {\Huge $\cdots$};
    \node at (15.5,1) {\Huge $\cdots$};

    \coordinate (ai) at (0,0);
    \coordinate (ai+1) at (2,0);
    \coordinate (ai+2) at (4,0);
    \coordinate (ai+3) at (6,0);
    \coordinate (ai+j-1) at (8,0);
    \coordinate (ai+j) at (10,0);
    \coordinate (ai+j+1) at (12,0);
    \coordinate (ai+j+2) at (14,0);

    \coordinate (bi) at (0,2);
    \coordinate (bi+1) at (2,2);
    \coordinate (bi+2) at (4,2);
    \coordinate (bi+3) at (6,2);
    \coordinate (bi+j-1) at (8,2);
    \coordinate (bi+j) at (10,2);
    \coordinate (bi+j+1) at (12,2);
    \coordinate (bi+j+2) at (14,2);    

    \node at (-1,3) (ci) {};
    \coordinate (ci+1) at (1,3);
    \coordinate (ci+2) at (3,3);
    \coordinate (ci+3) at (5,3);
    \coordinate (ci+4) at (7,3);
    \coordinate (ci+j) at (9,3);
    \coordinate (ci+j+1) at (11,3);
    \coordinate (ci+j+2) at (13,3);

        \draw[thick][gray!60, line width=1.7mm] (ci+3) .. controls (5,.5) and (7,1.3) .. (ai+j-1);
    \draw[thick] (ai) -- (ai+1) -- (ai+2) -- (ai+3);
    \draw[thick] (ai+j) -- (ai+j+2);
    \draw[thick] (ai) -- (bi);
    \draw[thick] (ai+1) -- (bi+1);
    \draw[thick] (ai+2) -- (bi+2);
    \draw[thick] (ai+3) -- (bi+3);
    \draw[thick][gray!60, line width=1.7mm] (ai+j-1) -- (bi+j-1);
    \draw[thick] (ai+j) -- (bi+j);
    \draw[thick] (ai+j+1) -- (bi+j+1);
    \draw[thick] (bi) -- (ci+1) -- (bi+1) -- (ci+2) -- (bi+2) -- (ci+3) -- (bi+3);
    \draw[thick][dashed, line width=.8mm] (bi+3) -- (ci+4);
    \draw[thick] (ci+4) -- (bi+j-1);
    \draw[thick] (bi+j-1) -- (ci+j) -- (bi+j) -- (ci+j+1) -- (bi+j+1) -- (ci+j+2) -- (bi+j+2) -- (ai+j+2);
    \draw[thick] (ci+1) .. controls (1,.5) and (3,1.3) .. (ai+2);
    \draw[thick] (ci+1) .. controls (1,.5) and (3,1.3) .. (ai+2);
    \draw[thick] (ci+2) .. controls (3,.5) and (5,1.3) .. (ai+3);

    \draw[thick] (ci+4) .. controls (7,.5) and (9,1.3) .. (ai+j);
    \draw[thick] (ci+j) .. controls (9,.5) and (11,1.3) .. (ai+j+1);
    \draw[thick] (ai) -- (bi) -- (ci+1) -- (bi+1) -- (ci+2);
    \draw[thick] (ai) -- (ai+3);
    \draw[gray!60, line width=1.7mm] (ai+3) -- (ai+j);
    \draw[thick] (ci+j) -- (bi+j) -- (ai+j);

      \draw[thick] (ci+j+1) .. controls (11,.5) and (13,1.3) .. (14,0);

    \node[fill=black, circle, inner sep=2pt] at (ai) {};
    \node[fill=black, circle, inner sep=2pt] at (ai+1) {};
    \node[fill=black, circle, inner sep=2pt] at (ai+2) {};
    \node[fill=black, circle, inner sep=2pt, label=-90:$v_2$] at (ai+3) {};
    \node[fill=black, circle, inner sep=2pt, label=-90:$a$] at (ai+j-1) {};
    \node[fill=black, circle, inner sep=2pt,
  label=-90:$v_4$] at (ai+j) {};
    \node[fill=black, circle, inner sep=2pt,
  label=-90:$v_5$] at (ai+j+1) {};
    \node[fill=black, circle, inner sep=2pt] at (ai+j+2) {};
    
    \node[fill=black, circle, inner sep=2pt] at (bi) {};
    \node[fill=black, circle, inner sep=2pt] at (bi+1) {};
    \node[fill=black, circle, inner sep=2pt] at (bi+2) {};
    \node[fill=black, circle, inner sep=2pt, label=-45:$x$] at (bi+3) {};
    \node[fill=black, circle, inner sep=2pt,
  label=-45:$v_3$] at (bi+j-1) {};
    \node[fill=black, circle, inner sep=2pt,
  label=-45:$v_7$] at (bi+j) {};
    \node[fill=black, circle, inner sep=2pt] at (bi+j+1) {};
    \node[fill=black, circle, inner sep=2pt] at (bi+j+2) {};
    
    \node[fill=black, circle, inner sep=2pt] at (ci+1) {};
    \node[fill=black, circle, inner sep=2pt] at (ci+2) {};
    \node[fill=black, circle, inner sep=2pt, label=90:$v_1$] at (ci+3) {};
    \node[fill=black, circle, draw=black, inner sep=2pt,
  label=90:$y$] at (ci+4) {};
    \node[fill=black, circle, inner sep=2pt,
  label=90:$v_6$] at (ci+j) {};
    \node[fill=black, circle, inner sep=2pt] at (ci+j+1) {};
    \node[fill=black, circle, inner sep=2pt] at (ci+j+2) {};    

      \begin{scope}
    \clip (-.5,3) rectangle (3.3,-.1);
    \draw[thick] (ci) .. controls (-1,.5) and (1,1.3) .. (2,0);
  \end{scope} 
  
  \begin{scope}
    \clip (12,3) rectangle (14.5,-.1);
    \draw[thick] (ci+j+2) .. controls (13,.5) and (15,1.3) .. (16,0);
  \end{scope}       
  \end{scope}

\node at (7.2,-7.2) {\Large $G_{m_1+m_2}$};
\end{tikzpicture}
}
    \caption{Isomorphism between $G_{m_1 + m_2}$ and graph obtained after subdividing the edge $xy$ (connecting two disconnected components $G_{m_1}$ and $G_{m_2}$). The dotted edge is being subdivided, and the gray edges are the new ones in the new graph.} 
    \label{fig:gm_subdivide_edge_different_component}
    \end{figure}

\begin{proof}
We know there is a bijection between the vertices of $H_{n-1}$ and $\Pp_n$; hence we just need to show that this bijection preserves edges. Let $\Ind(G)$ and $\Ind(G')$ be neighbors in $H_{n-1}$, and assume without loss of generality that $\Ind(G')$ is obtained by subdividing (instead of contracting) an edge $xy$ of $\Ind(G)$. Furthermore, assume $G = G_{m_1} \sqcup G_{m_2} \sqcup \dots \sqcup G_{m_s}$.  
By \cref{lem:edge_in_subdivision}, $x$ and $y$ must belong to two different components of $G$, say $G_{m_1}$ and $G_{m_2}$ without loss of generality. By \cref{l:noextraedges} and the fact that $G'$ is planar and ternary, we conclude that $G' \cong G_{m_1 + m_2} \sqcup G_{m_3} \sqcup \dots \sqcup G_{m_s}$. Therefore, the edge between $\Ind(G)$ and $\Ind(G')$ in $H_{n-1}$ gets mapped to the edge between $(m_1, \dots, m_s)$ and $(m_1 + m_2, m_3, \dots, m_s)$ in $\Pp_n$. 

Conversely, given the edge between $(m_1, \dots, m_s)$ and $(m_1 + m_2, m_3, \dots, m_s)$ in $\Pp_n$, we will show that there exists an edge between their corresponding vertices $\Ind(G_{m_1} \sqcup G_{m_2} \sqcup \dots \sqcup G_{m_s})$ and $\Ind(G_{m_1 + m_2} \sqcup G_{m_3} \sqcup \dots \sqcup G_{m_s})$ in $H_{n-1}$. Indeed, one can go from the first complex to the second by subdividing any $xy$ where $x \in G_{m_1}$ and $y \in G_{m_2}$ both having degree two. See \cref{fig:gm_subdivide_edge_different_component} for an illustration.
\end{proof}

\section{\texorpdfstring{The $f$, $h$ and $\gamma$-vectors}{The f, h and gamma-vectors} of planar flag spheres } \label{sec:VectorsPlanarFlagSpheres}

Recall that Gal's conjecture \cite{gal2005} posits that all flag homology spheres are $\gamma$-positive. Moreover, a real-rooted $h$-polynomial implies $\gamma$-positivity. Gal's conjecture was proposed after he found counterexamples to the stronger conjecture that the $h$-polynomials of flag spheres were real rooted. By \cite[Section 2.4]{gal2005}, flag simplicial spheres obtained by successive edge subdivisions from the crosspolytope are $\gamma$-positive. Therefore, \cref{thm:crosspolysubdivide} implies that $\Ind(G_m)$ is $\gamma$-positive. In this section, we investigate further the $f$-vector, $h$-vector and $\gamma$-vector of $\Ind(G_m)$, and prove the stronger result that a flag homology sphere whose $1$-skeleton is complementary to a planar graph has a real-rooted $h$-polynomial. 

For $m\ge 1$ and $0\le i\le m-1$, let $f_i(m) := f_i(\Ind(G_m))$. When needed, we use $f_m(m)=0$.

\begin{proposition}\label{p:face-pell-recursion}
For each $m\ge 1$, $f_0(m) = 3m-1$.
For each $m\ge 2$, $f_1(m) = (9m^2-19m+12)/2$.
For $m\ge 3$ and $2\le i\le m-1$,
$$
    f_i(m) = 2f_{i - 1}(m - 1) + f_i(m - 1) + f_{i - 2}(m - 2) + f_{i - 1}(m - 2).
$$
\end{proposition}
\begin{proof}
First note that $G_1=K_2$, and $\Ind(G_1)=2K_1$, so $f_0(1)=2$.  Also, $f_1(m)$ is the number of edges of $\Ind(G_m)$, that is, the number of pairs of vertices in $G_m$ that do not form an edge.  So for $m\ge 2$, $f_1(m) = \binom{3m-1}{2}-5(m-1) =(9m^2-19m+12)/2$. Now assume $m\ge 3$ and $2\le i\le m-1$.

    Refer to the recursive construction of $G_m$ for the vertices and edges of $G_m$ not in $G_{m-1}$.
 Note that $G_m - (N[a_m]\cup N[c_m]) = G_m - \{a_{m-1},b_{m-1},c_{m-1},a_m,b_m,c_m\} = G_{m - 2}$.

    Using the characterization of independent sets of $G_m$ from Proposition \ref{prop:indep_gm}, we count the number of independent sets of size $i+1$ for each of the five cases: Case $(1)$ gives exactly $f_i(m-1)$ subsets, Case $(2)$ gives $f_{i-1}(m-1)$ subsets, and Case $(3)$ gives exactly $f_{i-2}(m-2)$.
    
    For Cases~(4) and~(5), note that every independent set of size~$i$ of~$G_{m - 1}$ intersects 
$\{a_{m-1},\allowbreak b_{m-1},\allowbreak c_{m-1}\}$ in either 
$a_{m-1},\allowbreak b_{m-1},\allowbreak c_{m-1},\allowbreak \{a_{m-1},\allowbreak c_{m-1}\}$ 
or~$\emptyset$.
We now show that the number of these intersections is $f_{i-1}(m-1)+f_{i-1}(m-2)$. To see this, first consider independent sets $J$ of size $i$ (dimension $i-1$) in $G_{m-1}$.  If $J\cap\{a_{m-1},b_{m-1},c_{m-1}\}= \emptyset$, extend $J$ to a set of case 4 by adding $a_m$.  If $J\cap\{a_{m-1},b_{m-1},c_{m-1}\}\ne \emptyset$, then $J$ can be extended to a unique set of case 4 or 5 by adding $a_m$ (if $b_{m-1}\in J$) or $c_m$ (if $a_{m-1}$ and/or $c_{m-1}$ is in $J$).  This gives all independent sets of Case $(4)$ or $(5)$ except the sets $A$ for which $A\cap \{a_{m-1},b_{m-1},c_{m-1},a_m,b_m,c_m\} = \{c_m\}$.  The number of these is $f_{i-1}(m-2)$.
  
    Counting the elements in each partition we have 
    \begin{align*}
        f_i(m) & = f_i(m - 1) + f_{i - 1}(m - 1) + f_{i - 2}(m - 2) + f_{i - 1}(m - 1) + f_{i - 1}(m - 2) \\
               & = 2f_{i - 1}(m - 1) + f_i(m - 1) + f_{i - 2}(m - 2) + f_{i - 1}(m - 2).
    \end{align*}
\end{proof}

We now define the Delannoy polynomials mentioned in \cref{t:fhintro}.
\begin{definition}
The \em{Delannoy number} $D(m,k)$ is the number of lattice paths from $(0,0)$
to $(m,k)$ using steps $(1,0)$, $(0,1)$, and $(1,1)$.
\end{definition}
Following \cite{WZC2019},
let $d(m,k)=D(m-k,k)$, and for fixed $m$, let 
$d_m(t)$ be the generating function $d_m(t) = \sum_{k=0}^m d(m,k)t^k$.  It is known (see for example~\cite[A008288]{oeis} and~\cite[Section 2]{WZC2019}) that $d(m,k)$ satisfies the recursion $d(m,0) = d(m, m) = 1$ for $m \geq 0$ and, for $m > 1$ and $0 < k < m$, 
$$
    d(m,k) = d(m - 1,k) + d(m - 1, k - 1) + d(m - 2, k - 1).
$$
In 2019, Wang, Zheng and Chen~\cite{WZC2019} showed that Delannoy polynomials are real-rooted.

\begin{theorem}[{\cite[Theorem 2.1]{WZC2019}}]\label{thm:delannoyrealrooted}
The zeros of the Delannoy polynomial $\sum_{k=0}^m d(m,k)t^k$ are negative
real numbers.
\end{theorem}

We are now ready to tackle \cref{t:fhintro}. The first part of the theorem is stated for all planar graphs whose independence complexes are homology spheres. Let $G$ be such a graph and assume further it is connected. Then recall from \cref{rem:R3} that $G$ is either a $G_m$ or $R_3$. The case of $G = G_m$ is stated as the second part of \cref{t:fhintro} and proved as \cref{t:hvector-delannoy}. When $G = R_3$, the $h$-polynomial of $\Ind(G) = C_6$ is $t^2 + 4t+1$, which is real-rooted.

\begin{theorem}\label{t:hvector-delannoy}
Let $G$ be a planar ternary connected graph such that $\Ind(G)$ is an $(m-1)$-dimensional homology sphere. Then 
$h_k(\Ind(G)) = d(m,k)$.
\end{theorem}

\begin{proof}
We may assume $G = G_m$ by the previous discussion. Recall $f_i(m)= f_i(\Ind(G_m))$ and, similarly, write $h_k(m) = h_k(\Ind(G_m))$ and
$$h(m,t) = \sum_{k=0}^m h_k(m)t^{m-k} =
{\sum_{i=0}^m f_{i-1}(m)(t-1)^{m-i}}.$$
Proposition~\ref{p:face-pell-recursion}
gives
$$f_{i-1}(m) = 2f_{i - 2}(m - 1) + f_{i-1}(m - 1) + f_{i - 3}(m - 2) + f_{i - 2}(m - 2).$$
So
\begin{eqnarray*}
h(m,t) &=& \sum_{i=0}^m f_{i-1}(m)(t-1)^{m-i} \\
&=& 2\sum_{i=0}^m f_{i-2}(m-1)(t-1)^{m-i} +
 \sum_{i=0}^m f_{i-1}(m-1)(t-1)^{m-i} \\
 &\qquad& +\sum_{i=0}^m f_{i-3}(m-2)(t-1)^{m-i} +
 \sum_{i=0}^m f_{i-2}(m-2)(t-1)^{m-i}\\
&=& 2\sum_{i=0}^{m-1} f_{i-1}(m-1)(t-1)^{m-1-i} +
 (t-1)\sum_{i=0}^{m-1} f_{i-1}(m-1)(t-1)^{m-1-i} \\
 &\qquad& +\sum_{i=0}^{m-2} f_{i-1}(m-2)(t-1)^{m-2-i} +
 (t-1) \sum_{i=0}^{m-2} f_{i-1}(m-2)(t-1)^{m-2-i}\\
&=& (t+1)h(m-1,t)+th(m-2,t).
\end{eqnarray*}
Equating coefficients of $t^{m-k}$,
$h_k(m) = h_k(m-1)+h_{k-1}(m-1)+h_{k-1}(m-2)$.
This is the same recursion as for $d(m,k)$.
Since $(h_0(1),h_1(1)) = (1,1) = (d(1,0),d(1,1))$ and
$(h_0(2),h_1(2),h_2(2)) = (1,3,1) = (d(2,0),d(2,1),d(2,2))$,
$h_k(m) = d(m,k)$ for all $m\ge 1$, $0\le k\le m$.
\end{proof}

We immediately conclude the following from Theorems \ref{thm:delannoyrealrooted} and \ref{t:hvector-delannoy}.

\begin{corollary}
The zeros of the $h$-polynomial of $\Ind(G_m)$ are negative real numbers.
\end{corollary}

To remove the requirement of connectedness, recall from Sections \ref{sec:join} and \ref{sec:fhg} that for any pair of graphs $G$ and $H$,
\[
h(\Ind(G \sqcup H), t) = h(\Ind(G) \ast \Ind(H), t) = h(\Ind(G), t) \cdot h(\Ind(H), t).
\]
\begin{corollary}
If $G$ is a planar graph such that $\Ind(G)$ is a homology sphere, then the zeros of $h(\Ind(G), t)$ are negative real numbers.
\end{corollary}

\cref{t:fhintro} then follows directly from~\cref{t:hvector-delannoy} and the corollaries above.


\section{Future work}\label{sec:FutureWork}

Recall that Trung characterized Gorenstein planar graphs \cite{trung2018characterization}, all of which are ternary except the graph $R_3$. Pinter constructed Gorenstein nonplanar graphs in \cite{pinter1997}, but the ternary property of such graphs was not explored. In particular, many of the graphs introduced in~\cite{pinter1997} contain an induced subgraph isomorphic to $G_m$ for some $m$, see for example~\cite[Figure 4, Figure 10]{pinter1997}.
This naturally leads to the question of finding nonplanar ternary graphs that are Gorenstein. 

We now provide one such example below. Start with a disconnected graph $G^0= H_1 \sqcup H_2 \sqcup H_3$ with $H_i \cong C_5$ for $i \in [3]$. Consider vertices $1 \in V(H_1)$ and $6 \in V(H_2)$. We obtain a graph $G^1$ by subdividing the edge $\{1,6\}$ in $\Ind(G^0)$. Next, consider vertices $7 \in N_{H_2}(6)$ and $11 \in V(H_3)$. We obtain the graph $G^2$ pictured below by subdividing the edge $\{7,11\}$ in $\Ind(G^1)$.

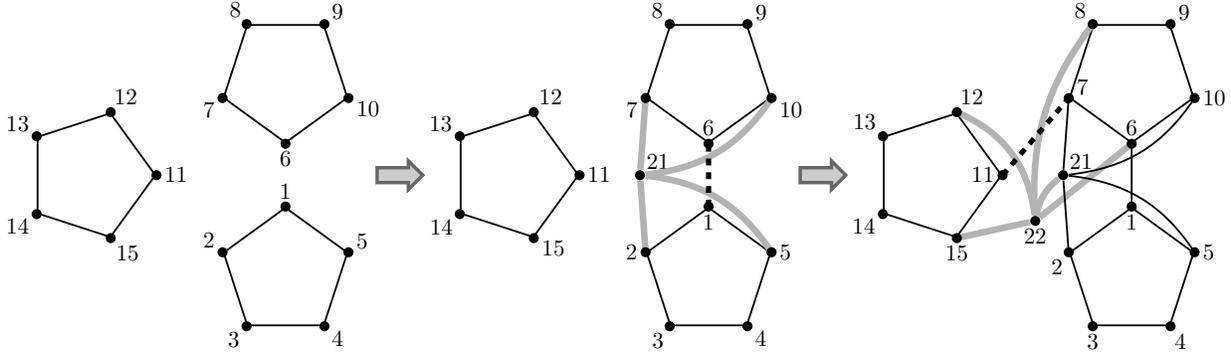
\begin{figure}[H]
    \centering
    \resizebox{1\textwidth}{!}{
\begin{tikzpicture}[font=\small]

\tikzset{
  Pent/.style={
    regular polygon,
    regular polygon sides=5,
    draw, thick,
    minimum size=2cm, inner sep=0pt,
    rotate=90         
  },
  vertex/.style={
    circle, fill=black, inner sep=0pt, minimum size=4.4pt
  },
  labstyle/.style={
    font=\footnotesize,
    inner sep=1pt
  },
  grayline/.style={
    draw=gray!60, line width=1mm
  }
}

  \begin{scope}[xshift=0cm]
    \node[Pent, rotate=180] (L1) at (0,0) {};
    \node[Pent, rotate=90]  (T1) at (3,1.5) {};
    \node[Pent, rotate=-90] (B1) at (3,-1.5) {};
    \foreach \i/\lab/\ang in {
      1/11/0, 2/12/80, 3/13/150, 4/14/210, 5/15/-60
    }{\node[vertex,label={[labstyle]\ang:\lab}] at (L1.corner \i) {};}
    \foreach \i/\lab/\ang in {
      1/6/-90, 2/10/-10, 3/9/30, 4/8/120, 5/7/200
    }{\node[vertex,label={[labstyle]\ang:\lab}] at (T1.corner \i) {};}
    \foreach \i/\lab/\ang in {
      1/1/90, 5/5/15, 4/4/-30, 3/3/220, 2/2/160
    }{\node[vertex,label={[labstyle]\ang:\lab}] at (B1.corner \i) {};}
  \end{scope}

 \node at (4.7,0) [single arrow,
                 draw=black!60,line width=1.5pt,
                 fill=gray!40,
                 minimum width=5pt,
                 single arrow head extend=3pt,
                 minimum height=7mm,inner sep=2.5pt,
                 rotate=0] {};

  \begin{scope}[xshift=6.5cm]
    \node[Pent, rotate=180] (L2) at (0,0) {};
    \node[Pent, rotate=90]  (T2) at (3,1.5) {};
    \node[Pent, rotate=-90] (B2) at (3,-1.5) {};

    \node[vertex, label={[labstyle]45:21}, xshift=-30, yshift=0] (c21) at ($(T2)!0.5!(B2)$) {};

    \draw[grayline]            (T2.corner 5) -- (c21) -- (B2.corner 2);
    \draw[grayline, bend left=25]  (T2.corner 2) to (c21);
    \draw[grayline, bend right=25] (B2.corner 5) to (c21);
    \draw[dashed, line width=.8mm] (B2.corner 1) to (T2.corner 1);

    \foreach \i/\lab/\ang in {
      1/11/0, 2/12/60, 3/13/150, 4/14/220, 5/15/-7
    }{\node[vertex,label={[labstyle]\ang:\lab}] at (L2.corner \i) {};}
    \foreach \i/\lab/\ang in {
      1/6/90, 2/10/-18, 3/9/15, 4/8/135, 5/7/200
    }{\node[vertex,label={[labstyle]\ang:\lab}] at (T2.corner \i) {};}
    \foreach \i/\lab/\ang in {
      1/1/-90, 2/2/180, 3/3/235, 4/4/330, 5/5/0
    }{\node[vertex,label={[labstyle]\ang:\lab}] at (B2.corner \i) {};}
  \end{scope}

\node at (11.2,0) [single arrow,
                 draw=black!60,line width=1.5pt,
                 fill=gray!40,
                 minimum width=5pt,
                 single arrow head extend=3pt,
                 minimum height=7mm,inner sep=2.5pt,
                 rotate=0] {};

  \begin{scope}[xshift=13cm]
    \node[Pent, rotate=180] (L3) at (0,0) {};
    \node[Pent, rotate=90]  (T3) at (3,1.5) {};
    \node[Pent, rotate=-90] (B3) at (3,-1.5) {};

    \node[vertex, label={[labstyle]45:21}, xshift=-30, yshift=0] (c21b) at ($(T3)!0.5!(B3)$) {};
    \node[vertex, inner sep=1pt, label={[labstyle]-90:22}, xshift=10, yshift=-20] (c22b) at ($(L3)!0.6!(c21b)$) {};

        \draw[grayline, bend left=30]  (L3.corner 2) to (c22b);
    \draw[grayline, bend right=30] (c22b) -- (L3.corner 5);
    \draw[grayline]                (T3.corner 1) to (c22b);
    \draw[grayline, bend right=20] (T3.corner 4) to (c22b);
    \draw[grayline, bend right=20] (c21b) to (c22b);
    \draw[dashed, line width=.8mm]                (L3.corner 1) to (T3.corner 5);
    
    \draw[thick]            (T3.corner 5) -- (c21b) -- (B3.corner 2);
    \draw[thick, bend left=25]  (T3.corner 2) to (c21b);
    \draw[thick, bend right=25] (B3.corner 5) to (c21b);
    \draw[thick] (B3.corner 1) to (T3.corner 1);

    \foreach \i/\lab/\ang in {
      1/11/180, 2/12/80, 3/13/135, 4/14/210, 5/15/-90
    }{\node[vertex,label={[labstyle]\ang:\lab}] at (L3.corner \i) {};}
    \foreach \i/\lab/\ang in {
      1/6/90, 2/10/0, 3/9/20, 4/8/135, 5/7/10
    }{\node[vertex,label={[labstyle]\ang:\lab}] at (T3.corner \i) {};}
    \foreach \i/\lab/\ang in {
      1/1/-90, 2/2/235, 3/3/270, 4/4/-75, 5/5/0
    }{\node[vertex,label={[labstyle]\ang:\lab}] at (B3.corner \i) {};}
  \end{scope}

\end{tikzpicture}
}
    \caption{Constructing a nonplanar ternary Gorenstein graph}
    \label{fig:ex_nonplanar_Gorenstein}
\end{figure}
It is quick to check that all induced cycles in $G^2$ have length four or five. $G^2$ is nonplanar as the induced subgraph of $G^2$ on vertices $\{6, 7, 8, 9, 10, 21, 22\}$ is homeomorphic to $K_{3,3}$. Lastly, to see that $G^2$ is a Gorenstein graph, notice that $\Ind(G^0)$ is homeomorphic to $\Ind(G^2)$ since we obtain $\Ind(G^2)$ from $\Ind(G^0)$ by subdividing edges, and by~\cref{thm:ternarygorenstein} $\Ind(G^0)$ is homotopy equivalent to $S^{\dim \Ind(G^0)}$. The result then follows since $\dim \Ind(G^0) = \dim \Ind(G^2)$.

We can use a generalization of this construction to obtain Gorenstein graphs, and we suspect that many of these are ternary. Additionally, under certain conditions we can guarantee the output of this construction to be nonplanar. The general construction looks as follows.

\begin{enumerate}
    \item Set $k=0$.
    \item Let $G^k \cong \bigsqcup_{i=1}^n H_i$ with $H_i \cong G_{m_i}$ for all $i$ and $n \geq 2$.
    \item Pick two vertices $v_i \in V(H_i)$ and $v_j \in V(H_j)$ such that $v_i, v_j$ are in different components of $G^k$.
    \item Subdivide the edge $v_iv_j$ in $\Ind(G^k)$ and call the resulting graph $G^{k+1}$.
    \item If $G^{k+1}$ is connected, stop the procedure. Else, either stop the procedure or set $k=k+1$ and return to (3).
\end{enumerate}

Observe that we used a special case of this in the proof of Theorem \ref{thm:partitionrefine}. In particular, as long as we choose the vertices $v_i$ and $v_j$ such that $\deg v_i = \deg v_j = 2$, the graph $G^k$ is isomorphic to a disjoint union of $G_m$, that is, it is a planar ternary Gorenstein graph.

We can actually show that $G^k$ is nonplanar if and only if we at least once pick some vertex $v_i$ such that $\deg v_i \geq 3$: Note that the forward direction is given by the above observation. For the other direction, let $k$ be the largest such that at each prior step $\deg v_i = \deg v_j = 2$. Hence, we currently have $G^k \cong \bigsqcup_{i=1}^{n} H_i$ with $H_i \cong G_{m_i}$ for all $i$. Now assume $\deg v_i \geq 3$ and let $N(v_i) \supseteq \{x_1, x_2, x_3\}$. Calling the new vertex in $G^{k+1}$ $y$, the induced subgraph on vertices $\{y, x_1, x_2, x_3, v_i\}$ is isomorphic to $K_{2,3}$. By the structure of $G_m$ we can additionally find a vertex $z \in V(H_i)$ such that for each $\ell \in [3]$ there exist pairwise disjoint paths from $x_\ell$ to $z$ in $H_i$ avoiding $v_i$. Hence, $G^{k+1}$ contains a subgraph homeomorphic to $K_{3,3}$ and is therefore nonplanar.

\begin{question} 
        If $G^k$ is nonplanar, under which conditions is it also ternary?
\end{question}

\begin{remark}
    Suppose $G^0$ has $N$ connected components, and the construction above produces a graph $G$.  Form a graph $W$ on $N$ vertices, each corresponding to a connected component of $G^0$, with two vertices adjacent if at some point the construction uses vertices from the corresponding components.
    We believe that if $W$ is a tree, then the graph $G$ is ternary.
\end{remark}

\begin{question}
    How do the conditions for being nonplanar or ternary change if we just require $v_i$ and $v_j$ to be in different components of $G^k$, but not necessarily vertices of $H_i$?
\end{question}

\subsection*{Acknowledgements}

This project started during the Graduate Research Workshop in Combinatorics 2024, which was hosted by the University of Wisconsin--Milwaukee and supported in part by NSF Grant DMS-195344. The authors would like to thank Aram Bingham, Jeremy Martin and Shira Zerbib for their early contributions to this project. The authors also thank the anonymous referee for several helpful suggestions. Richard Danner's research was partially supported by NSF Grant DMS-2246967 and Simons Collaboration Gift \#854037.  Yirong Yang's research was partially supported by a Graduate Fellowship from NSF Grant DMS-2246399.

\end{document}